\documentclass[twocolumn]{autart}    % Enable this line and disable the 
\usepackage{subcaption}
\usepackage{amsmath} 
\usepackage{amssymb}  % assumes amsmath 
\usepackage{mathtools}
\usepackage{bbold}
\usepackage{algorithm}
\usepackage{algpseudocode}
\usepackage{comment}
\usepackage{xcolor}
\usepackage{graphicx,epstopdf}

\begin{document}

\begin{frontmatter}

\title{On Linear Quadratic Potential Games\thanksref{footnoteinfo}} % Title, preferably not more 

\thanks[footnoteinfo]{This paper appears on Automatica \\ {https://doi.org/10.1016/j.automatica.2025.112643} }

\author[UBC]{Sara Hosseinirad}
% \ead{sarahrad@ece.ubc.ca},    % Add the 
\author[EPFL]{Giulio Salizzoni}
% \ead{giulio.salizzoni@epfl.ch},
\author[UBC]{Alireza Alian Porzani}
% \ead{yaphetsf75g@gmail.com},               % e-mail address 
  % (ead) as shown
\author[EPFL]{Maryam Kamgarpour}
% \ead{maryam.kamgarpour@epfl.ch}  % (ead) as shown

\address[UBC]{Department of Electrical and Computer Engineering, The University of British Columbia}  % Please supply                                           
\address[EPFL]{SYCAMORE Lab, School of Engineering, EPFL}             % full addresses
          
                             % keyword list or with the 
                                          % help of the Automatica 
                                          % keyword wizard

\begin{abstract}                          % Abstract of not more than 200 words.
Our paper addresses characterizing conditions for linear quadratic (LQ) games to be potential. The desired properties of potential games in finite action settings, such as convergence of learning dynamics to Nash equilibria and challenges of learning Nash equilibria in continuous state and action settings, motivate us to characterize LQ potential games. Our first contribution is to consider two-player LQ games with full-state feedback and scalar states and action space, and to analytically verify that the set of potential games within this example is limited, essentially differing only slightly from an identical interest game. Given this finding, we restrict the class of LQ games to those with decoupled dynamics and decoupled state linear feedback information structure. For this subclass, we show that the set of potential games strictly includes non-identical interest games and characterize conditions for LQ games in this subclass to be potential. We further derive their corresponding potential function and prove the existence of a Nash equilibrium. Meanwhile, we highlight the challenges in characterizations of the Nash equilibrium for this class of potential LQ games, theoretically and through simulation.
\end{abstract}

\end{frontmatter}

\section{Introduction}

\label{sec:introduction}
% general setup
Numerous emerging networked dynamic systems, such as smart grids and autonomous vehicles, rely on the decision-making of multiple interacting agents. Due to the increasing complexities of the dynamics and uncertainties of these systems, data-driven control has gained increasing attention. Reinforcement learning approaches like policy gradient methods have been widely applied in single-agent settings, and their convergence properties have been theoretically investigated~\cite{Fazel18}. However, multi-agent reinforcement learning has been less understood in both theory and practice.  

% game theoretic
Non-cooperative \textit{linear quadratic} (LQ) games serve as a fundamental setting for understanding the convergence of multi-agent reinforcement learning algorithms. This class of games admits Nash equilibria in linear policies, which can be computed by solving a set of coupled Riccati equations~\cite{Basar98}. Given these well-understood properties of LQ games, several recent studies have sought to theoretically understand multi-agent reinforcement learning algorithms within the LQ game setting. 

% A multi-agent control problem can be formulated as a non-cooperative game in which each agent (also referred to as a player) aims at optimizing its loss function. 
% This function depends on the state of the system and the actions of the other agents; hence, there is a coupling between the agent's decisions.
Most prior works on learning in LQ games have focused on policy gradient approaches. While the convergence of policy gradient and its variants in single-agent \textit{infinite-horizon} LQ Gaussian control problems is well-established~\cite{Fazel18}, the case of LQ games has proven more challenging. In particular, in infinite-horizon general-sum LQ game settings, ~\cite{Mazumdar20a} and~\cite{Mazumdar20b} argued that there are neither global nor local convergence guarantees to a Nash equilibrium for policy gradient. This negative result has motivated several follow-up studies on understanding the convergence of policy gradient in LQ games. 

For \textit{finite-horizon} general-sum LQ games, the global convergence of the natural policy gradient method to the Nash equilibrium is proven in~\cite{Hambly21}. However, this result relies on the assumption of sufficient noise for exploration. 
One study showed that the policy gradient algorithm has convergence guarantees when agents interact through the aggregated state and action~\cite{Roudneshin20}. Other studies considered zero-sum \cite{zhang2019policy}, mean-field \cite{carmona2019linear}, or risk-sensitive \cite{zhao2023global} LQ games, proving the convergence of policy gradient variants in these settings. 
% can not be easily extended to cases where there is a restriction on information structure, i.e., not all the agents have full knowledge of all the states. The convergence of policy gradient methods for these LQ games remains an open question. 

In \textit{Markov} games, of which LQ games are a subset, the challenges in proving convergence to Nash equilibria have motivated researchers to consider specific subclasses, particularly potential dynamic games. In static games, potential games are known for several desirable properties such as the existence of pure strategy Nash equilibria, the convergence of the best response dynamics~\cite{monderer1996potential}, and learning dynamics to such equilibria \cite{heliou2017learning}. These advances also have solid practical relevance in real-world dynamic games. The LQ game is a well-known dynamic game setting for many real-world applications, including congestion games \cite{rosenthal1973class}, resource allocation games~\cite{Lambart2018}, games arising in multi-vehicle autonomous driving~\cite{Aghajani2015,Jiang20}, and electricity markets~\cite{chen22,Paccagnan2016}. Consequently, a natural approach to extend the convergence guarantees of learning-based methods to dynamic games is to characterize the set of dynamic games that admit a potential function.

% These advances have also solid practical relevance in that several classes of static games encountered in the real world do admit a potential function. 
% ~\cite{Jiang20}
% electricity markets CITE:Pang/Kanzow...  
% MK finite state and actions Markov potential games: 1) what kind of games are potential; 2) what has been proven for PG in these games
%MK: note: "Global Convergence of Localized Policy Iteration in Networked Multi-Agent Reinforcement Learning" focuses on single-agent objective but localized policies. they manage to prove convergence by assuming agents communicate .. 
%MK: https://proceedings.mlr.press/v216/zhou23b/zhou23b.pdf
% network potential games... perhaps can be extended to LQG
Recent literature in the multi-agent learning community has focused on conditions that ensure a dynamic game with finite state and action spaces is potential. Studies on Markov games with finite states and actions~\cite{leonardos22,Zhang23,Zhang22} reveal that deriving conditions to characterize a potential Markov game is not straightforward. In particular,~\cite{leonardos22} demonstrated that the stage cost being potential is neither necessary nor sufficient for the dynamic game to be potential. Given these limitations,~\cite{Zhang23} focused on Markov games with decoupled dynamics and local policies where agents make decisions using only their local state. In this setting, they proved the convergence of a distributed learning algorithm to a Nash equilibrium for Markov potential games. By characterizing LQ potential games as a subset of Markov potential games, we gain a fundamental understanding of a significant class of games in control theory.

As a first step toward understanding Markov potential games in infinite state and action spaces, we consider LQ games. Existing results on LQ potential games are limited. In~\cite{Zazo16}, a theoretical analysis of dynamic potential games with constrained state-action sets for LQ games is provided by extending the definition of static potential games. In~\cite{Aprem19}, Bayesian optimization is used to develop a novel algorithm for computing the Nash equilibrium in LQ potential games. Both studies provided some examples of LQ potential games beyond identical interest games. However, in our current work, we challenge their characterization of potential games (see Remark~\ref{rem:pastLQpotential}) aiming to clarify a misunderstanding in previous works on the characterization of potential games.  
% There were a few works for infinite state and action we looked at.. we should briefly summarize what they showed. 

% We should highlight the following from [10] to be consistent with the line of logic: deriving conditions for a Markov game to be potential is difficult. In particular, stage-wise... Given this challenge, [20] focused on xxx. then, end this paragraph by saying to our knowledge, the conditions for potential LQ games were not derived in past work. Then, start a new paragraph about what we do in this paper...
% The setting is similar to the one covered in this study, with a distinction where the LQ games have continuous action and state spaces. The convergence of policy gradient algorithms in Markov potential games has been investigated in several studies~\cite{Zhang23},~\cite{leonardos22},~\cite{Zazo18}, and ~\cite{Zhang22}. Here, the convergence of these algorithms in the LQ game potential game settings is studied. 

% but in their case, the sets of states and actions are finite, while these are not finite for the LQ games as states and actions are continuous. 
% Moreover, the characteristic of this type of game is an important discussion that has been overlooked in these studies, and we will cover it extensively in this work, particularly in the context of LQ games.
\subsection{Contributions}
Generally, the objective of this paper is to advance a fundamental understanding of classes of LQ potential games. Our contributions towards this goal are as follows:
\begin{itemize}
    \item We derive conditions under which finite-horizon general-sum LQ games admit a potential function (see Lemma~\ref{Lm:Lemma7}). Leveraging this result,  we show that the LQ potential games with coupled dynamics and a full-state feedback information structure form a limiting class, meaning they are "almost" the same as identical interest games (see Proposition \ref{prop:coupled_i} for precise statement). 
    \item Motivated by the aforementioned limitations, we consider a subclass of LQ games with agents having decoupled dynamics and decoupled information structure but coupled loss functions. For this subclass, we characterize the set of LQ potential games and show that this set is significantly larger than the set of identical interest games (see Theorem~\ref{Th:CL_po}). 
    \item For potential games within this subclass, we derive the potential function (Proposition~\ref{prop:pot_function}) and show the existence of Nash equilibria (Proposition~\ref{prop:NEminimum}). Furthermore, we formalize the challenges associated with proving the uniqueness of the Nash equilibrium and its computation theoretically (Proposition~\ref{prop:initialStatesCorrelation} and Section~\ref{sec:challenge}) and through simulation (Section~\ref{sec:simulation}). 
    
%\item Given limitations in computing Nash equilibria, we provide conditions under which a policy gradient algorithm can provably converge to a stationary point of the potential function. 
%    \item We provide conditions under which a policy gradient algorithm can provably converge to a stationary point of the potential function.
\end{itemize} 

The rest of the paper is organized as follows. In Section~\ref{sec:preliminaries}, we introduce the class of games under study and provide background information on potential games. In Section~\ref{sec:potential}, we derive conditions for LQ potential games with coupled dynamics and full-state feedback and characterize the potential games in the case of decoupled dynamics and information structure. The properties of the LQ potential game, its corresponding structured optimal control problem, and the convergence of policy gradient for this control problem to a stationary point are provided in Section~\ref{decoupled_game}. We conclude and discuss open research directions in Section~\ref{sec:Conclusion}.

% In Section~\ref{sec:gradDescProof}, a model-based policy gradient algorithm for the potential games with decoupled dynamics and the decoupled information structure is proposed, and we prove that this algorithm converges to a stationary point in the class of potential games. 

% Through simulation, we verify the convergence of this algorithm to a stationary point. 
\subsection{Notations}
Herein, $\mathbb{R}$, $\mathbb{R}_{\geq 0}$ and $\mathbb{R}_{> 0}$ refer to the sets of real numbers,  non-negative, and positive real numbers, respectively. The summation of natural numbers $\{n^p, \dots, n^q\}$ is denoted by $n^{p:q} = \sum_{j = p}^{q}n^j$. The set $[N]$ represents \{1, \dots, N\}, and $[N] \backslash \{i\}$ denotes all members of the set except $i$. The zero vector is denoted by $\mathbb{0}_n \in \mathbb{R}^n$, and the one vector by $\mathbb{1}_n \in \mathbb{R}^n$. We write $A = \text{blockdiag}\{A^1, \dots, A^N\}$ to denote a block-diagonal matrix with $A^1, \dots, A^N$ on its diagonal entries. 
% If we are given matrices $A \in \mathbb{R}^{m \times n}$, $B \in \mathbb{R}^{p \times q}$, and $C \in \mathbb{R}^{m \times n}$, the Kronecker product is represented by $(A \otimes B) \in \mathbb{R}^{mp \times nq}$, and the Hadamard product by $(A \odot C) \in \mathbb{R}^{m \times n}$.
For a vector $(\gamma^1, \dots, \gamma^N) \in \mathbb{R}^{n}$ where $n = n^1 + \dots + n^N$, the vector $\gamma^{-i}\in \mathbb{R}^{\sum_{j \neq i }^{N}n^j}$ is defined as $\gamma^{-i}: = (\gamma^1, \dots, \gamma^{i-1},\gamma^{i+1}, \dots, \gamma^N) $. The identity matrix of dimension $n$ is represented by $I_n$. Given a matrix $Q\in \mathbb{R}^{n \times n}$, $(Q)_{ij}$ is a sub-matrix of $Q$ that is formed by selecting entries from rows $(n^{1:i-1} +1)$ to $n^{1:i}$ and columns $(n^{1:j-1} +1)$ to $n^{1:j}$. 
% The matrix $E^i_n \in \mathbb{R}^{n\times n}$ is defined such that its entries are zero, except for the $i$th column where $E_{ji} = 1$ for all $j \in \{1, \dots, n\}$. The matrix $D^i_n \in \mathbb{R}^{n\times n}$ is defined such that its entries are zero except for $D_{ii} = 1$, which is the $i$th diagonal one.
 % Given $v$, $\{v^i_{\tau}\}_{\tau=0}^{t} = (v^i_{0}, \dots, v^i_{t})$. 

\section{Linear quadratic game setup}
\label{sec:preliminaries}
We define a finite-horizon {linear quadratic} (LQ) game and a variant of it with agents having decoupled dynamics. The notion of Nash equilibrium under different information structures and the definition of LQ potential games are provided. Our first result, Lemma \ref{Lm:Lemma7}, derives conditions for an LQ game to be potential. This lemma serves as a foundation for deriving necessary and sufficient conditions for LQ potential games throughout this work. 
\subsection{Linear quadratic games}
Consider a non-cooperative general-sum LQ game where the dynamics are defined as
\begin{equation} \label{eq:accum_dynamic}
	x_{t+1} = Ax_{t} + \sum_{i = 1}^{N}B^iu^i_{t}, \quad x_0 \sim \mathcal{D},
\end{equation}
with $t \in \{0, \dots, t_f-1\}$. For simplicity, let us consider scalar actions denoted by $u^i_t \in \mathbb{R}$, and the states of the game are $x_t \in \mathbb{R}^{n}$. The sub-index and super-index indicate the time and the agent, respectively. The initial state $x_0 \in \mathbb{R}^{n}$ is randomly distributed according to a distribution $\mathcal{D}$. We assume that all agents choose their decision variables simultaneously at each time step.

The action of agent $i$, $u^i_t$, depends on the information available to the agent. Let us denote the decision variable of agent $i$ as a function of its information by $\gamma^i_t \in \mathbb{R}^{q^i}$. The exact form of $\gamma^i$ will be detailed in Section~\ref{sec:info_struct}. The loss function of each agent is
\begin{flalign} \label{eq:cost_lqg}
	J^i(\gamma) =& \mathop{\mathbb{E}}_{x_0 \sim \mathcal{D}} \bigl[J^i_x(\gamma,x_0) + J^i_u(\gamma,x_0)\bigr],
 % \bigg[(x_{t_f}-d_{t_f})^T Q_{t_f}^i(x_{t_f}-d_{t_f}) +\\ \nonumber
 % &   \sum_{t = 0}^{t_f-1} \bigg((x_t-d_{t})^T Q_t^i(x_t-d_t) + \\  
 % & \big(\sum_{j = 1}^{N}u_t^j R_t^{ji}+ M_t\big)u^i_t\bigg)\bigg]
\end{flalign}
where the state part $J^i_x(\gamma,x_0) \in \mathbb{R}$ and the action part $J^i_u(\gamma,x_0) \in \mathbb{R}$ of the loss function are defined as follows:
\begin{flalign} \label{eq:Loss_state_gen}
    J^i_x(\gamma,x_0) :=& \sum_{t = 0}^{t_f} (x_t-d_{t})^T Q_t^i(x_t-d_t), \\ \label{eq:Loss_action_gen}
    J^i_u(\gamma,x_0) :=& \sum_{t = 0}^{t_f-1} \sum_{j = 1}^{N}\sum_{h = 1}^{N} (u_t^j)^T (R^i_t)_{jh}u^h_t,
\end{flalign}
with $\gamma^i = (\gamma^i_0, \dots, \gamma^i_{t_f-1})\in \mathbb{R}^{t_fq^i}$ and $\gamma = (\gamma^1, \dots, \gamma^N)\in \mathbb{R}^{t_fq^{1:N}}$. The desired state at time $t$ is denoted by $d_t \in \mathbb{R}^{n}$.  The parameters $Q_t^i \in \mathbb{R}^{n \times n}$ for $t \in \{0, \dots, t_f\}$ are positive semi-definite. For $t \in \{0, \dots, t_f-1\}$ and for $i,j,h \in [N]$, the cross-term action loss coefficients are $(R^i_t)_{jh} = (R^i_t)_{hj} \in \mathbb{R}$, and other action loss coefficients are $(R^i_t)_{ll} \in \mathbb{R}_{\geq 0}$ for  $l \in [N] \backslash \{i\}$ and $(R^i_t)_{ii} \in \mathbb{R}_{> 0}$.

% where $u^i = (u^i_0, \dots, u^i_{t_f-1})$, $u = (u^1, \dots, u^N)$, and 
% and the stage loss is $\forall t \in \{0, \dots, t_f-1\}$
% \begin{flalign}\label{eq:cost_ct_lqg}
% 	c^i_t(x_t, u^1_t, \dots, u^N_t)  = &   \big[x_t^T Q_t^ix_t + (\sum_{j = 1}^{N}u_t^j R_t^{ji}+ d_t)u^i_t\big], \\
%     c^i_{t_f}(x_{t_f})  = &\mathop{\mathbb{E}}_{x_0 \sim \mathcal{D}} \big[x_{t_f}^T Q_{t_f}^ix_{t_f} \big].
% \end{flalign}

% For the derivation of potential functions, the game is assumed to be deterministic, i.e., the noise distribution and the random initial state are avoided in the game definition. 
We say that the LQ game in~(\ref{eq:accum_dynamic}) and~(\ref{eq:cost_lqg}) has \emph{decoupled dynamics} if for $i \in [N]$ and $\ t \in \{0, \dots, t_f-1\}$,
\begin{flalign}\label{eq:dyn_decoupled}
    x_{t+1}^i =& A^ix^i_{t} + b^iu_{t}^i, 
\end{flalign}
where $x_t^i \in \mathbb{R}^{n^i}$, $A^i \in \mathbb{R}^{n^i \times n^i}$, $b^i \in \mathbb{R}^{n^i}$ with $\sum_{i =1}^{N}n^i = n$.
Note that the decoupled dynamics can be written as the dynamics introduced in~(\ref{eq:accum_dynamic}) by defining the joint state as $x_t = (x_t^1, \dots, x_t^N)$, joint control as $u_t = (u_t^1, \dots, u_t^N)$,  and  $A = \text{blockdiag}\{A^1, \dots, A^N\} \in \mathbb{R}^{n \times n}$, $B^i = (\mathbb{0}_{n^{1:i-1}}, b^i, \mathbb{0}_{n^{i + 1:N}}) \in \mathbb{R}^{n}$,  for $i \in [N]$. 

The motivation for considering LQ games with decoupled dynamics is to derive less restrictive conditions for the existence of a potential function. As demonstrated in the following two examples, this class of LQ games also holds practical relevance.

% The evolution of each agent's state can be solved independently not knowing other players' states; however, the loss function of a player depends on the state of others. 
% Note that an LQ game with decoupled dynamic is not an identical-interest game as matrices $Q_t^i$ and $R_t^{ij}$ are not the same for all players.
% \begin{exmp}
\textbf{Example 1}:
In the context of decentralized formation control of multi-vehicle systems, each vehicle can be considered an agent, with $x_t^i$ and $u_t^i$ representing the position and control action of the $i$th vehicle at time step $t$. As noted in~\cite{Aghajani2015}, one approach to defining each vehicle's objective is to minimize the formation error and energy consumed by the vehicle. In~\cite{Aghajani2015}, the formation error of the $i$th vehicle is defined as follows:
\begin{align*}
e^i_x(x^1_t, \dots, x^N_t) \ = \sum_{j = 1}^{N} w^{ij}_t \| x_t^i - x^j_t - d^{ij}_t \|^2,
\end{align*}
where $d_t^{ij} = d^i_t - d^j_t$ indicates the desired distance between two vehicles. The error weights are $w^{ij}_t \in \mathbb{R}_{\geq 0}$. The energy consumption of vehicle $i$ at time $t$ is presented by the quadratic form of control actions $(u^i_t)^T(R^i_t)_{ii}u^i_t$ where $(R^i_t)_{jh}$ is nonzero if and only if  $j = i$ and $h = i$. The agents' objective is the summation of formation error and energy consumption over a finite horizon and thus, can be written as in~(\ref{eq:cost_lqg}).
% \end{exmp}
% by defining a positive semi-definite matrix $Q^i_t \in \mathbb{R}^{n \times n}$ as
% \begin{equation*}
%     Q^i_t = ({W}_t^i +W^i_tE^i_N + {E^i_N}^TW^i_t + \sum_{j = 1}^{N} w_t^{ij} D^i_N)\otimes I_{n/N},
% \end{equation*}
% where $W_t^i = \text{blockdiag} \{w^{i1}, \dots, w^{iN}\} \in \mathbb{R}^{N \times N}$.

% \begin{example}
% \label{exp:cournotGame}
\textbf{Example 2}: Let us consider a simple dynamic Cournot game, where agent $i$ is a firm that decides how much product $u^i_t \in \mathbb{R}$ to sell at each time step, given the amount of product stored $x^i_t \in \mathbb{R}_{\geq0}$. A simple abstraction of the agent's dynamics is
% At each time step, the agent also receives a certain amount of product $a^i_t \in \mathbb{R}$, which however is not a control variable:
\begin{align*}
    x^i_{t+1} = x^i_t - u^i_t.
\end{align*}
%where $0\leq x^i_t \leq x_{\textbf{max}}$, and $x_{\textbf{max}}$ is the maximum products that can be stored.
The firm's goal is to maximize its profit which is a function of the price of the product. This price is a decreasing function of the total production of all agents. As commonly done in a Cournot model in~\cite{Shi20},~\cite{Lambart2018}, and~\cite{Paccagnan2016}, let us consider a linear price at time $t$ of $p_t(u_t) = -\alpha_t \sum_{i=1}^N u^i_t$
% \begin{align*}
%     p_t(u_t) = -\alpha_t \sum_{i=1}^N u^i_t
%     % + p_{0,t},
% \end{align*}
where $\alpha_t \in \mathbb{R}_{>0}$. Maximizing the profit is equivalent to minimizing the following loss function:
\begin{align}\label{eq:cournot_loss}
    J^i(\gamma) =& \mathop{\mathbb{E}}_{x_0 \sim \mathcal{D}} \bigl[ Q_{t_f}^i(x^i_{t_f}-d_{t_f}^i)^2 - \sum_{t=0}^{t_f-1}p_t(u_t)u^i_t\bigr].
\end{align}
 The loss function parameters $Q_t^i = 0$ for  $t\in \{0, \dots, t_f-1\}$. In~(\ref{eq:Loss_state_gen}), the final stage loss parameter is $Q_{t_f}^i \in \mathbb{R}_{>0}$, and $d_{t_f}$ is the desired stored products at the final stage. In~(\ref{eq:Loss_action_gen}), the action loss coefficients $(R^i_t)_{ji} = \alpha_t$ are the same for  $i, j \in [N]$, and the rest of the action loss coefficients are $(R^i_t)_{jh} = 0$ for  $j,h \in [N] \backslash \{i\}$. 
 % Commonly, there is a linear term in the loss function of the Cournot game, where we consider it as zero without the loss of generality.
 % There is also a linear term $-p_{0,t}u^i_t$.
% \end{example}

\subsection{Nash equilibria and potential game definitions} \label{sec:info_struct}
For the LQ games defined in~(\ref{eq:accum_dynamic}) and~(\ref{eq:cost_lqg}) we consider three different information structures: (i) open-loop, (ii) full-state linear feedback, and (iii) decoupled state
linear feedback. In each case, for  $i \in [N]$, the decision variables are\\
\textbf{(i)} actions $\gamma_t^i = u_t^i \in \mathbb{R}$ with $q^i = 1$;\\
\textbf{(ii)} full-state linear feedback coefficients $\gamma_t^i = K^i_t \in \mathbb{R}^{n}$, resulting in $u_t^i = -K_t^ix_t$ and $q^i =n$;\\
\textbf{(iii)} decoupled state linear feedback coefficients $\gamma_t^i = k^i_t \in \mathbb{R}^{n^i}$, resulting in $u_t^i = -k_t^ix^i_t$ and $q^i = n^i$.

% \begin{itemize}
%     \item (i) 
%     \item (ii) 
%     \item (iii) 
% \end{itemize}
%   \smallskip
\begin{defn}\label{def:ne}
    The joint strategy $({\gamma}^{1*}, \dots, {\gamma}^{N*}) \in \mathbb{R}^{t_fq^{1:N}}$ is a Nash equilibrium (NE) for the LQ game defined in~(\ref{eq:accum_dynamic}) and~(\ref{eq:cost_lqg}) if and only if for $i \in [N]$,
\begin{equation*}
% \label{eq:NE_CL}
		J^i(\gamma^{i*},\gamma^{-i*}) \leq J^i(\tilde{\gamma}^{i},\gamma^{-i*}), \qquad \forall \tilde{\gamma}^{i}\in \mathbb{R}^{t_fq^i}.
\end{equation*}
\end{defn} 
% \smallskip
\begin{defn}\label{def:potential}
  The LQ game defined in~(\ref{eq:accum_dynamic}) and~(\ref{eq:cost_lqg}) is a potential game if and only if there exists a function $\Pi: \mathbb{R}^{t_fq^{1:N}} \rightarrow \mathbb{R}$ such that  for $i \in [N]$, $ \gamma^{-i} \in \mathbb{R}^{ t_f\sum_{j \neq i}^{N}q^j}$, and $\hat{\gamma}^{i}, {\gamma}^i \in \mathbb{R}^{ t_fq^i}$, 
	\begin{eqnarray} \label{eq:PotentialDef}
		J^i(\gamma^i, \gamma^{-i}) - J^i( \hat{\gamma}^i,\gamma^{-i}) = 
		\Pi(\gamma^i, \gamma^{-i}) - \Pi(\hat{\gamma}^i,\gamma^{-i}).
	\end{eqnarray} 
In this case, $\Pi$ is called a potential function for the game. 
\end{defn}
% Although any potential function maximizer is a Nash equilibrium, a potential game might have Nash equilibria that are not maximizers of the potential function.
Identical interest LQ games defined below are trivial examples of potential games. While their definition is apparent from their name, to clarify in comparing our results we formalize them below.
\begin{defn}\label{def:identical_interest}
    The LQ games in~(\ref{eq:accum_dynamic}) and~(\ref{eq:cost_lqg}) are identical interest if for $i,j,h,l \in [N]$ and for $t \in \{0, \dots, t_f - 1\}$,
\begin{flalign} \label{eq:ID_QR}
    Q^i_{t+1} = Q^j_{t+1}, \quad \text{and} \quad
    (R^i_t)_{lh} = (R^j_t)_{lh}.
\end{flalign}
\end{defn}
% With an example, we demonstrated that the condition in~(\ref{eq:ID_QR}) is not necessary for the LQ game to be an identical interest game. 
% Assume an LQ game as described in~(\ref{eq:accum_dynamic}) and~(\ref{eq:cost_lqg}) with full state information structure where $t_f = 1$, $N = 2$. The decision variable of this game is $\gamma = (K_0^1, K_0^2)$. The loss function of agent 1 can be written as
% \begin{flalign*}
%     &J^1(\gamma, x_0) = \mathop{\mathbb{E}}_{x_0 \sim \mathcal{D}} \bigl[x_0^TQ_0^1x_0 +\\
%    & (Ax_0 + B^1u_0^1 + B^2u_0^2)^TQ_1^1(Ax_0 + B^1u_0^1 + B^2u_0^2)\\
%    & + (R_0^1)_{11}(u_0^1)^2 + (R_0^1)_{12}u_0^1u_0^2 + (R_0^1)_{22}(u_0^2)^2 \big]
% \end{flalign*}
% where $x_1$ is replaced by  $Ax_0 - B^1u_0^1 - B^2u_0^2$. This game is identical interest if and only if $J^1(\gamma, x_0) = J^2(\gamma, x_0)$. This equality is met if and only if 
% \begin{flalign}
    
% \end{flalign}
Our focus is on deriving conditions to identify non-trivial LQ potential games. Therefore, we derive the necessary and sufficient conditions for an LQ game to be potential.

\section{LQ potential games} \label{sec:potential}
First, in Lemma~\ref{Lm:Lemma7}, we develop a test to verify whether an LQ game with the three information structures above is potential. Using this lemma,  we demonstrate that, unlike the open-loop information structure in (i), a game with full-state feedback information structure (ii) is potential if and only if agents' losses are identical at each stage, except for the initial stage. This motivates us to consider decoupled dynamics and information structure (iii) defined in Subsection~\ref{sec:potential_decoupled}. For this subset of LQ games, we characterize conditions under which the game admits a potential function, beyond the identical interest loss function. 

\subsection{Necessary and sufficient conditions for LQ potential games} 
\label{sec:pre}
%Our goal is to leverage this result to define necessary and sufficient conditions for LQ games defined in the previous section to be differentiable.

The game \emph{pseudo-gradient}, $\mathcal{G}(\gamma): \mathbb{R}^{t_fq^{1:N}} \to \mathbb{R}^{t_fq^{1:N}}$, is defined as the gradient of each agent loss function with respect to her own decision variable as follows:
\begin{align*}
% \label{eq:pseudogradient}
    \mathcal{G}(\gamma) = \begin{bmatrix}
        \big(\frac{\partial J^1 (\gamma)}{\partial \gamma^1}\big)^T, \dots, 
        \big(\frac{\partial J^N (\gamma)}{\partial \gamma^N}\big)^T
    \end{bmatrix}^T,
\end{align*}
where $\frac{\partial J^i (\gamma)}{\partial \gamma^i} \in \mathbb{R}^{t_fq^i}$. For a differentiable $\mathcal{G}(\gamma)$, the \emph{Jacobian} of the pseudo-gradient, $\mathcal{J}: \mathbb{R}^{t_fq^{1:N}} \to \mathbb{R}^{t_fq^{1:N} \times t_fq^{1:N}}$, is:
\begin{align} \label{eq:jac_j}
    \mathcal{J}(\gamma) = \begin{bmatrix}
		\frac{\partial^2 J^1(\gamma)}{\partial \gamma^1 \partial \gamma^1} & \dots & \frac{\partial^2 J^1(\gamma)}{\partial \gamma^N \partial \gamma^1}  \\
                \vdots & \ddots & \vdots \\
		\frac{\partial^2 J^N(\gamma)}{\partial \gamma^1 \partial \gamma^N} & \dots&  \frac{\partial^2 J^N(\gamma)}{\partial \gamma^N \partial \gamma^N}
	\end{bmatrix},
\end{align}
where $\frac{\partial^2 J^i(\gamma)}{\partial \gamma^j\partial \gamma^i} \in \mathbb{R}^{t_fq^i \times t_fq^j}$.
Our first result develops an equivalent condition to an LQ game being potential.
% \smallskip
\begin{lem} \label{Lm:Lemma7}
The LQ game defined in~(\ref{eq:accum_dynamic}) and~(\ref{eq:cost_lqg}) is a potential game if and only if for ${\gamma}_t^i \in \mathbb{R}^{q^i}$, ${\gamma}_{\tau}^j \in \mathbb{R}^{q^j}$, $i,j\in [N]$, and $t, \tau \in \{0, \dots, t_f-1 \}$
\begin{equation} \label{eq:Lemma7}
	\frac{\partial^2 J^i(\gamma)}{\partial \gamma_{\tau}^j\partial \gamma_t^i} =  \frac{\partial^2 J^j(\gamma)}{\partial \gamma_{\tau}^j \partial \gamma_t^i}.
\end{equation}
 % \smallskip
\end{lem}
\begin{pf}
A static differentiable game is potential if and only if the Jacobian of its pseudo-gradient is symmetric~\cite[Theorem 1.3.1]{Facchinei03}.  Our goal is now to verify that the condition in \eqref{eq:Lemma7} is equivalent to the game Jacobian being symmetric. The block $ij$ of the matrix in~(\ref{eq:jac_j}) is:
    \begin{equation*}
    % \label{eq:subJac}
	\frac{\partial^2 J^i(\gamma)}{\partial \gamma^j\partial \gamma^i} = 
	\begin{bmatrix}
		\frac{\partial^2 J^i(\gamma)}{\partial \gamma_{0}^j\partial \gamma_0^i} & \dots & \frac{\partial^2 J^i(\gamma)}{\partial \gamma_{t_f-1}^j\partial \gamma_0^i}  \\
                \vdots & \ddots & \vdots \\
		\frac{\partial^2 J^i(\gamma)}{\partial \gamma_{0}^i\partial \gamma_{t_f-1}^i} & \dots&  \frac{\partial^2 J^i(\gamma)}{\partial \gamma_{t_f-1}^j\partial \gamma_{t_f-1}^i}
	\end{bmatrix}.
    \end{equation*} 
    Since $J^i(\gamma)$ is a polynomial function of $\gamma$, it is twice-continuously differentiable, and thus
    \begin{flalign*}
        \frac{\partial^2 J^i(\gamma)}{\partial \gamma^j \partial \gamma^i} = \left[ \frac{\partial^2 J^i(\gamma)}{\partial \gamma^i \partial \gamma^j}\right]^T.
    \end{flalign*}
    Hence, the terms on the diagonal of the Jacobian in~(\ref{eq:jac_j}) are symmetric. Since diagonal terms are symmetric, the Jacobian is symmetric if and only if
    \begin{flalign*}
        \frac{\partial^2 J^i(\gamma)}{\partial \gamma^j \partial \gamma^i} = \left[ \frac{\partial^2 J^j(\gamma)}{\partial \gamma^i \partial \gamma^j}\right]^T.
    \end{flalign*}
Note that the left-hand side of equation~(\ref{eq:Lemma7}) is the entry in the $t$-th row and $\tau$-th column of $\frac{\partial^2 J^i(\gamma)}{\partial \gamma^j \partial \gamma^i}$ and the right-hand side is the entry in the $t$-th row and $\tau$-th column of $\frac{\partial^2 J^j(\gamma)}{\partial \gamma^i \partial \gamma^j}$, and $\frac{\partial^2 J^j(\gamma)}{\partial \gamma^i \partial \gamma^j}$ is symmetric. Thus, the equality above holds if and only if conditions in~(\ref{eq:Lemma7}) are met. \qed
    \end{pf}
 % entries of these two matrices shown in~(\ref{eq:subJac}) are equals. Therefore, the Jacobian is symmetric if and only if the condition in equation~(\ref{eq:Lemma7}) is met.   
% Investigating condition~(\ref{eq:Lemma7}) can be challenging, especially for games with full state feedback information structure (ii) and where $u_t^i \in \mathbb{R}^{n_u}$ for all agents. This is because computing the first derivative of the loss function with respect to $K_{t}^i \in \mathbb{R}^{n_u \times n}$ results in a $n_{u} \times n$ matrix. In order to obtain the left-hand side of~(\ref{eq:Lemma7}), the second derivative of the loss function needs to be computed with respect to a $n_{u} \times n$ matrix, $K_{t}^j$, which can be a challenging task.
%\begin{lemma} \label{Lm:commmon}
%The general-sum finite-horizon LQ game defined in~(\ref{eq:accum_dynamic}) and~(\ref{eq:cost_lqg}) is potential if and only if the loss function of each agent $i \in [N]$ can be expressed as the sum of a term that is identical to all agents (identical interest game) plus another term that does not depend on its own policy (dummy game):

Note that each element of $\mathcal{J}(\gamma)$ is a second derivative of the loss function $J^i(\gamma)$, and the gradients of the loss functions can be more than second-order polynomials with respect to the decision variable, depending on the information structure. Therefore, investigating conditions for desirable properties in the Jacobian, such as diagonally strictly concave~\cite{rosen65}, can be as challenging as deriving LQ potential conditions. Next, we focus on investigating the potential condition in Lemma~\ref{Lm:Lemma7} for different game settings and information structures.

\subsection{Potential games with coupled dynamics} 
\label{sec:pot_nondecoupled}

Lemma~\ref{Lm:Lemma7} provides a test for verifying whether an LQ is potential. However, computing the needed gradients in state-feedback information structure in general is challenging. Thus, we focus on a scalar state setting. In this setting, for full-state feedback information structure, we show that the class of potential games is almost the same as the class of identical interest games (see conditions in Proposition \ref{prop:coupled_i}).

%Assume an example of the LQ game with a non-decoupled dynamic and scalar state.
%Through a simple example, we develop necessary and sufficient conditions under which the game described in the example with an open-loop information structure (i) is potential. 
%Moreover, we establish that the necessary and sufficient condition for the game with open-loop information structure in (i) to be potential is solely a necessary condition for the game with information structure in (ii) to be potential.
% \subsubsection*{Example} 
\textbf{Example 3.} Here we consider an LQ game defined by equations~(\ref{eq:accum_dynamic}) and~(\ref{eq:cost_lqg}) where $N=2$, $n = 1$, and $d_t^1 = d_t^2 = 0$, and the dynamics is $x_{t+1} = ax_t + b^1u^1_t + b^2u^2_t$ for  $t\in \{0, \dots, t_f-1\}$. In loss function in~(\ref{eq:cost_lqg}), for open-loop information structure (i), $\gamma_t^i= u_t^i$, and for full-state feedback (ii), $\gamma_t^i= K_t^i$ where $u_t^i = -K_t^ix_t$, for $i \in \{1, 2\}$. 
\begin{prop} \label{prop:coupled_i}
    For the setting described in Example 3,
    \begin{itemize}
        \item[(a)] the game is an open-loop potential game (information structure (i)) if and only if
        
          \begin{align} \label{eq:prop5}
        Q^1_1b^1b^2 + (R^1_0)_{12} = Q^2_1b^1b^2 + (R^2_0)_{12}.
\end{align}
and for $t \in \{1, \dots, t_f-1\} $
\begin{flalign} \label{eq:prop5_a}
            Q^1_{t+1} = Q^2_{t+1}, \quad \text{and} \quad
            (R_t^1)_{12} = (R_t^2)_{12}.
        \end{flalign}

        \item[(b)] the game is potential with full-state feedback information structure (ii) if and only if for $t=0$, equation~(\ref{eq:prop5}) holds, and for $t \in \{1, \dots, t_f-1\}$ and for $l,h \in \{1, 2\}$,
        \begin{flalign} \label{eq:prop5_b}
            Q^1_{t+1} = Q^2_{t+1}, \quad \text{and} \quad
            (R^1_t)_{lh} = (R^2_t)_{lh}.
        \end{flalign}
    \end{itemize}
 
\end{prop}
\begin{rem} \label{rem:pastLQpotential} Let us clarify a misunderstanding in the past works on LQ potential games. Past works~\cite{Aprem19,Zazo16} claimed that the LQ games in~(\ref{eq:cost_lqg}) with $Q_t^i = Q_t^j$, $(R^i_t)_{jh} = 0$ for $j\neq i$ and $h \neq i$, and any $(R^i_t)_{ii}$ for  $i \in [N]$, are LQ potential game under full-state feedback structure. However, our Proposition~\ref{prop:coupled_i} demonstrates that such an LQ game is an open-loop potential but not potential with a full-state feedback structure. Specifically, the loss parameters of quadratic terms also need to be identical in the latter case, i.e., $(R^i_t)_{ii} = (R^j_t)_{ii}$ for $i,j \in [N]$. This condition is not met in their provided game setting since $(R^j_t)_{ii} = 0$ for $j\neq i$ and $(R^i_t)_{ii}$ is not necessarily zero.
\end{rem}

Before presenting the proof, some insight into the derived results is offered. Notably, for $t>0$, the stage losses must be identical for the information structure (ii) whereas for the information structure (i) the loss of agent $i$ on its control action need not be the same as that of another agent. The condition on stage loss parameters for the first step $t=0$ is different than for $t> 0 $ due to the varying dependency of stage losses on $\gamma_t$, for $t >0$. In particular, the first stage loss only depends on $\gamma_0$, and the last stage loss depends on $\gamma_t$, for $t \in \{0, \dots, t_f-1\}$. Consequently, the conditions on the parameters of the first stage loss appear only when evaluating the condition in~(\ref{eq:Lemma7}) for $(t, \tau) = (0, 0)$.

%  In particular, contrary to the intuitive view that identical interest LQ games are potential, our proposition shows that this is not the case. 
% This point has also been misunderstood in past literature.

\begin{pf}
For scalar states and actions, the loss function of agent $i$ in~(\ref{eq:cost_lqg}) where $i \in \{1,2\}$ for $t_f = 2$ is 
\begin{flalign} \nonumber
    J^i(\gamma) &= \mathop{\mathbb{E}}_{x_0 \sim \mathcal{D}} \Big[ \sum_{t = 0}^{2} Q_t^i(x_t)^2 + \sum_{t=0}^{1} \Big((R_t^i)_{11}(u_t^1)^2\\
    &  + 2(R_t^i)_{12}u_t^1u_t^2 + (R_t^i)_{22}(u_t^2)^2 \Big) \Big], \label{eq:loss_scalar}
\end{flalign}
where $x_{t+1} = a^tx_0 + \sum_{\tau=0}^{t}a^{t-\tau}[b^1u^1_{\tau} + b^2u^2_{\tau}]$ for open-loop information structure (i), and $x_{t+1} = \prod_{\tau = 0}^{t}(a - b^1K_{\tau}^1 - b^2K^2_{\tau})x_0$ for full-state feedback (ii). By leveraging  Lemma \ref{Lm:Lemma7}, the above decomposition is used in the proof. 

\textit{(a) Open-loop information structure (i)}: First,
we show that for $t_f = 2$, the game is potential if and only if~(\ref{eq:prop5}) is met for $t = 0$, and conditions in~(\ref{eq:prop5_a}) are satisfied for $t = 1$. Subsequently, through induction, we demonstrate the validity of the proposition for any finite $t_f > 2$.

For $t_f = 2$, we verify that the condition in~(\ref{eq:Lemma7}) is satisfied for the four pairs of time steps $(t,\tau)$ where $t, \tau \in \{0, 1\}$, for $i = 1$, and for $j = 2$ if and only if conditions in Proposition~\ref{prop:coupled_i} part (a) is held. In this proof, the derivatives of states are utilized, $\frac{\partial x_{1}}{\partial u_0^i} = b^i$,  $\frac{\partial x_{1}}{\partial u_1^i} = 0$, $\frac{\partial x_{2}}{\partial u_0^i} = ab^i$, and $\frac{\partial x_{2}}{\partial u_1^i} = b^i$  for  $i \in \{1,2\}$. Observe that
\begin{flalign*}
    \frac{\partial J^1(\gamma)}{\partial u_1^1} &= 2\mathop{\mathbb{E}}_{x_0 \sim \mathcal{D}} \left[ Q_{2}^1x_{2} \frac{\partial x_{2}}{\partial u_1^1} + (R^1_1)_{11}u^1_1+ (R^1_1)_{12}u^{2}_1 \right].
\end{flalign*}
The second derivatives of $J^1(\gamma)$ are
\begin{flalign*}
    \frac{\partial^2 J^1(\gamma)}{\partial u_1^2 \partial u_1^1} &= 2 \bigl(Q_{2}^1b^1b^2 + (R_1^1)_{12}\bigr),\\
    \frac{\partial^2 J^1(\gamma)}{\partial u_0^2 \partial u_1^1} &= 2Q_{2}^1ab^1b^2.
\end{flalign*}
If $(R_1^1)_{12}$ and $Q_{2}^1$ are respectively replaced by $(R_1^2)_{12}$ and $Q_{2}^2$ in the equations above, the second derivatives of $J^2(\gamma)$ with respect to $u^1_1$ and $u^2_{\tau}$, for  $\tau \in \{0, 1\}$ are achieved. Condition in~(\ref{eq:Lemma7}) is met for two pairs of $(t,\tau) \in \{(1, 1), (1, 0)\}$ if and only if $Q_{2}^1ab^1b^2 = Q_{2}^2ab^1b^2$, and $Q_{2}^1b^1b^2 + (R_1^1)_{12} = Q_{2}^2b^1b^2 + (R_1^2)_{12}$. Since states are scalar, and dynamics parameters $b^1$ and $b^2$ are nonzero, these two conditions are met if and only if $Q_2^1 = Q_2^2$ and $(R_1^1)_{12} = (R_1^2)_{12}$ that are equivalent to~(\ref{eq:prop5_a}) for $t = 1$. 

% The derivative of $J^1(\gamma)$ and $J^2(\gamma)$ with respect to $u_0^1$ is
% \begin{flalign*}
%     \frac{\partial J^1(\gamma)}{\partial u_0^1} &= 2\mathop{\mathbb{E}}_{x_0 \sim \mathcal{D}} \bigl[Q_{2}^1x_{2} \frac{\partial x_{2}}{\partial u_0^1} + Q_{1}^1x_{1} \frac{\partial x_{1}}{\partial u_0^1} \\
%     & + (R^1_0)_{11}u^1_0+ (R^1_0)_{12}u^{2}_0 \big]. \\
%     \frac{\partial J^2(\gamma)}{\partial u_0^1} &= 2\mathop{\mathbb{E}}_{x_0 \sim \mathcal{D}} \bigl[ Q_{2}^2x_{2} \frac{\partial x_{2}}{\partial u_0^1} + Q_{1}^2x_{1} \frac{\partial x_{1}}{\partial u_0^1} \\
%     & + (R^2_0)_{11}u^1_0+ (R^2_0)_{12}u^{2}_0 \big].
% \end{flalign*}
% where $\frac{\partial x_{2}}{\partial u_0^i} = ab^i$ and $\frac{\partial x_{1}}{\partial u_0^i} = b^i$ for all $i \in \{1,2\}$.

The derivative of $J^1(\gamma)$ with respect to $u_0^1$ is
\begin{flalign*}
    \frac{\partial J^1(\gamma)}{\partial u_0^1} &= 2\mathop{\mathbb{E}}_{x_0 \sim \mathcal{D}} \biggl[Q_{2}^1x_{2} \frac{\partial x_{2}}{\partial u_0^1} + Q_{1}^1x_{1} \frac{\partial x_{1}}{\partial u_0^1} \\
    & + (R^1_0)_{11}u^1_0+ (R^1_0)_{12}u^{2}_0 \biggr].
\end{flalign*}
One can simply check that for the pair of  $(t, \tau) = (0, 1)$ in~(\ref{eq:Lemma7}), the repetitive condition $Q_2^1 = Q_2^2$ will be derived. The second derivative of $J^1(\gamma)$ with respect to $u_0^2$ is 
\begin{flalign*}
    \frac{\partial^2 J^1(\gamma)}{\partial u_0^2 \partial u_0^1} &= 2 \bigl[(a)^2Q^1_2b^1b^2 + Q^1_1b^1b^2 + (R^1_0)_{12} \big].
\end{flalign*}
If $(R_0^1)_{12}$, $Q_{2}^1$, and $Q_{1}^1$ are respectively replaced by $(R_0^2)_{12}$, $Q_{2}^2$, and $Q_{1}^2$ in the equation above, the second derivatives of $J^2(\gamma)$ with respect to $u^1_0$ and $u^2_0$ is achieved.
Since we already established that $Q_2^1 = Q_2^2$, for $(t, \tau) = (0,0)$, the condition in~(\ref{eq:Lemma7}) is satisfied if and only if 
\begin{flalign*}
    Q^1_1b^1b^2 + (R^1_0)_{12} = Q^2_1b^1b^2 + (R^2_0)_{12},
\end{flalign*}
where is identical to~(\ref{eq:prop5}).
By consolidating all necessary and sufficient conditions for different pairs of $(t, \tau)$, the Proposition~\ref{prop:coupled_i} part (a) is proved for $t_f =2$.

The proof can be extended for $t_f>2$ through induction. We assume that conditions in Proposition~\ref{prop:coupled_i} part (a) are necessary and sufficient conditions for the LQ game in Example 2 with $t_f = h$ with open-loop information structure to be potential, and we aim to demonstrate that these conditions are also necessary and sufficient conditions for the game with $t_f = h+1$. The difference between the loss function of agent 1 for $t_f = h$ and for $t_f = h+1$ lies in the last stage loss $Q_{h+1}^1(x_{h+1})^2 + (R_h^1)_{11}(u_h^1)^2 + 2(R_h^1)_{12}u_h^1u_h^2 + (R_h^1)_{22}(u_h^2)^2$. The decision variables $u_h^1$ and $u_h^2$ only appear in the last stage loss so that the second derivative of $J^1(\gamma)$ is
\begin{flalign*}
    \frac{\partial^2 J^1(\gamma)}{\partial u_h^2 \partial u_h^1} &= 2 \big[Q_{h+1}^1b^1b^2 + (R_h^1)_{12}\big],\\
    \frac{\partial^2 J^1(\gamma)}{\partial u_{h-1}^2 \partial u_h^1} &= 2Q_{h+1}^1ab^1b^2.
\end{flalign*}
Therefore, for the pair of $(t, \tau) \in \{(h, h), (h, h-1)\}$, the condition in~(\ref{eq:Lemma7}) is met if and only if $Q_{h+1}^1 = Q_{h+1}^2$ and $(R_h^1)_{12} = (R_h^2)_{12}$. The contribution of the last stage loss in the conditions obtained by other pairs of $(t, \tau)$ is repetitive, and other stage losses result in similar conditions derived for the case where $t_f = h$. Therefore, the game described in Example 3 with $t_f =h+1$ is open-loop potential if and only if conditions in Proposition~\ref{prop:coupled_i} part (a) are satisfied.

\textit{(b) Full-state feedback information structure (ii)}:
First, we prove that for $t_f = 2$, the game is potential if and only if equation~(\ref{eq:prop5}) is met for $t = 0$, and conditions in~(\ref{eq:prop5_b}) are satisfied for $t = 1$. Then we utilize induction to demonstrate that Proposition~\ref{prop:coupled_i} part (b) is also true for $t_f > 2$.

For $t_f = 2$, we examine that the condition in~(\ref{eq:Lemma7}) is satisfied for the four pairs of time steps $(t,\tau)$ where $t, \tau \in \{0, 1\}$  for $i = 1$, and for $j = 2$ if and only if conditions in Proposition~\ref{prop:coupled_i} part (b) are held.
In this proof, the derivatives of states, $\frac{\partial x_{1}}{\partial K_0^i} = -b^ix_0$,  $\frac{\partial x_{1}}{\partial K_1^i} = 0$, $\frac{\partial x_{2}}{\partial K_0^i} =  -b^i(a- b^1K_1^1 - b^2K_1^2)x_0$, and $\frac{\partial x_{2}}{\partial K_1^i} = -b^ix_1$, and nonzero derivatives of actions, $\frac{\partial u^i_{1}}{\partial K_0^j} = K^i_1b^jx_0$,  $\frac{\partial u^i_{1}}{\partial K_1^i} = -x_1$, and $\frac{\partial u^i_{0}}{\partial K_0^i} = -x_0$,
are utilized for  $i, j \in \{1,2\}$. Here, we only derive $J^1(\gamma)$ derivatives, but similar to the open-loop case, the $J^2(\gamma)$ derivative can be computed by replacing the super-index 1 with 2 for $Q$ and $R$ parameters. The derivative of $J^1(\gamma)$ in~(\ref{eq:loss_scalar}) with respect to $K_1^1$ is
\begin{flalign*}
    \frac{\partial J^1(\gamma)}{\partial K_1^1} &= 2\mathop{\mathbb{E}}_{x_0 \sim \mathcal{D}} \biggl[Q_{2}^1x_{2} \frac{\partial x_{2}}{\partial K_1^1} \\
    & + \{(R^1_1)_{11}u^1_1+ (R^1_1)_{12}u^{2}_1 \}\frac{\partial u^1_{1}}{\partial K_1^1} \biggr].
\end{flalign*}
The derivative of $\frac{\partial J^1(\gamma)}{\partial K_1^1}$ with respect to $K_1^2$ and $K_0^2$ are
\begin{flalign*}
    \frac{\partial^2 J^1(\gamma)}{\partial K_1^2 \partial K_1^1} &= 2 [Q_{2}^1b^1b^2 + (R_1^1)_{12}]\mathop{\mathbb{E}}_{x_0 \sim \mathcal{D}} \left[(x_1)^2\right],\\
    \frac{\partial^2 J^1(\gamma)}{\partial K_0^2 \partial K_1^1} &= 4 \big[Q_{2}^1b^1b^2(a -b^1K_1^1 - b^1K_1^2)\\
    &- b^2\{(R_1^1)_{11}K_1^1 + (R_1^1)_{12}K_1^2\}\big]\mathop{\mathbb{E}}_{x_0 \sim \mathcal{D}} [x_1x_0].
\end{flalign*}
For $(t, \tau) = (1, 0)$, the condition in~(\ref{eq:Lemma7}) is  met if and only if the coefficient of $K_1^1$ and $K_1^2$ and the constant term of $J^1(\gamma)$ and $J^2(\gamma)$ second derivatives are equal as follows.
\begin{flalign*}
    Q_2^1(b^1)^2b^2 &= Q_2^2(b^1)^2b^2\\
    Q_2^1(b^1)^2b^2 + b^2(R_1^1)_{11} &= Q_2^2(b^1)^2b^2 + b^2(R_1^2)_{11}\\
    Q_2^1b^1(b^2)^2 + b^2(R_1^1)_{12} &= Q_2^2b^1(b^2)^2 + b^2(R_1^2)_{12}.
\end{flalign*}
Thus, the condition in~(\ref{eq:Lemma7}) is met for $(t, \tau) = (1, 0)$ if and only if $Q_{2}^1 = Q_{2}^2$, $(R_1^1)_{12} = (R_1^2)_{12}$, and $(R_1^1)_{11} = (R_1^2)_{11}$, which are comparable with~(\ref{eq:prop5_b}) for $t = 1$. Consequently, the second derivatives of $J^1(\gamma)$ and $J^2(\gamma)$ with respect to $K_1^1$ and $K_1^2$ are also equal. 

The first derivative of $J^1(\gamma)$ in~(\ref{eq:loss_scalar}) with respect to $K_0^1$ is
\begin{flalign*}
    \frac{\partial J^1(\gamma)}{\partial K_0^1} &= 2\mathop{\mathbb{E}}_{x_0 \sim \mathcal{D}} \biggl[Q_{2}^1x_{2} \frac{\partial x_{2}}{\partial K_0^1} + Q_{1}^1x_{1} \frac{\partial x_{1}}{\partial K_0^1} \\
    & + \{(R^1_0)_{11}u^1_0 + (R^1_0)_{12}u^{2}_0 \}\frac{\partial u^1_{0}}{\partial K_0^1} \\
    & +\{(R^1_1)_{11}u^1_1 + (R^1_1)_{12}u^{2}_1 \}\frac{\partial u^1_{1}}{\partial K_0^1} \\ & + \{(R^1_1)_{12}u^{1}_1 + (R^1_1)_{22}u^2_1 \} \frac{\partial u^2_{1}}{\partial K_0^1}\biggr].
\end{flalign*}
The derivative of $\frac{\partial J^1(\gamma)}{\partial K_0^1}$ with respect to $K_1^2$ and $K_0^2$ are
\begin{flalign*}
    \frac{\partial^2 J^1(\gamma)}{\partial K_1^2 \partial K_0^1} &= 4 \big[Q_{2}^1b^1b^2(a - b^1K_1^1 - b^2K_1^2) \\
    &  - b^1\{(R_1^1)_{12}K_1^1 + (R_1^1)_{22}K^2_1\}\big]\mathop{\mathbb{E}}_{x_0 \sim \mathcal{D}} [x_1x_0]\\
    \frac{\partial^2 J^1(\gamma)}{\partial K_0^2 \partial K_0^1} &= 2 \{b^1b^2\big[Q_{1}^1 + Q_2^1(a - b^1K_1^1 - b^2K_1^2)^2 \\
    & + (R_1^1)_{11}(K_1^1)^2 + 2(R_1^1)_{12}K^1_1K^2_1 \\
    &+ (R_1^1)_{22}(K_1^2)^2\big]+ (R_0^1)_{12}\}\mathop{\mathbb{E}}_{x_0 \sim \mathcal{D}} [(x_0)^2].
\end{flalign*}
Since we showed that $Q_{2}^1 = Q_{2}^2$, $(R_1^1)_{12} = (R_1^2)_{12}$, and $(R_1^1)_{11} = (R_1^2)_{11}$, the condition in~(\ref{eq:Lemma7}) is met for $(t, \tau) = (0, 1)$ if and only if $(R_1^1)_{22} = (R_1^2)_{22}$. The equality in~(\ref{eq:Lemma7}) is also true for $(t, \tau) = (0, 0)$ if and only if
\begin{flalign*}
    Q^1_1b^1b^2 + (R^1_0)_{12} = Q^2_1b^1b^2 + (R^2_0)_{12},
\end{flalign*}
which is equivalent to Equation~(\ref{eq:prop5}). Consolidating conditions derived for all pairs of $(t, \tau)$, we notice the game is potential with full-state feedback information structure if and only if equation~(\ref{eq:prop5}) is met for $t = 0$, and conditions in~(\ref{eq:prop5_b}) are met for $t \in \{1,2\}$.

The proof can be extended for $t_f>2$ using induction.
We assume that conditions in Proposition~\ref{prop:coupled_i} part (b) are necessary and sufficient conditions for the LQ game in Example 2 with $t_f = h$ with full-state feedback information structure to be potential. Then, we aim to demonstrate that these conditions are necessary and sufficient for the game with $t_f = h+1$. The difference between $J^1(\gamma)$ for $t_f = h$ and $t_f = h+1$ is the last stage loss $Q_{h+1}^1(x_{h+1})^2 + (R_h^1)_{11}(u_h^1)^2 + 2(R_h^1)_{12}u_h^1u_h^2 + (R_h^1)_{22}(u_h^2)^2$. One can simply check that 
\begin{flalign*}
    \frac{\partial^2 J^1(\gamma)}{\partial K_{h-1}^2 \partial K_h^1} &= 4\big[Q_{h+1}^1ab^1b^2 \\
    & - \{Q_{h+1}^1(b^1)^2b^2 + b^2(R_h^1)_{11}\}K_h^1 \\
    & - \{Q_{h+1}^1b^1(b^2)^2 + b^2(R_h^1)_{12}\}K_{h-1}^2\big]\\
    & \times \mathop{\mathbb{E}}_{x_0 \sim \mathcal{D}} [x_{h-1}x_{h-2} ],\\
    \frac{\partial^2 J^1(\gamma)}{\partial K_h^2 \partial K_{h-1}^1} &= 4 \big[Q_{h+1}^1b^1b^2(a - b^1K_h^1 - b^2K_h^2) - b^1 \\
    & \{(R_{h}^1)_{12}K_h^1+ (R_h^1)_{22}K^2_h\}\big]\mathop{\mathbb{E}}_{x_0 \sim \mathcal{D}}[x_hx_{h-1}].
\end{flalign*}
Therefore, for the pair of $(t, \tau) = (h, h-1)\}$, the condition in~(\ref{eq:Lemma7}) is met if and only if $Q_{h+1}^1 = Q_{h+1}^2$, $(R_h^1)_{11} = (R_h^2)_{11}$, $(R_h^1)_{22} = (R_h^2)_{22}$, and $(R_h^1)_{12} = (R_h^2)_{12}$.
\qed
\end{pf}
%The following remark indicates the difference between the conditions presented in Proposition~\ref{prop:coupled_i} part (a) and part (b).

%In an open-loop information structure, the loss function is a second-order polynomial function of decision variables. On the other hand, with the full-state feedback information structure, the
%loss function depends on the decision variables with a polynomial of degrees higher than two.

\begin{comment}
\begin{rem}
In contrast to common and perhaps intuitive view that identical interest LQ games are potential, our proposition shows that this is not the case. In particular, in studies ~\cite{Zazo16} and~\cite{Aprem19}, the LQ games with $Q_t^i = Q_t^j$ and $(R^i_t)_{ij} = 0$ are provided as examples of potential LQ game with full-state feedback structure. However, our Proposition~\ref{prop:coupled_i},  demonstrates that such an LQ game is an open-loop potential and is not potential with a full-state feedback structure. 
\end{rem}
\end{comment}
%Example 3 represents a simple class of LQ games with scalar states and input variables. Despite its simplicity, it demonstrated that very restrictive conditions should be generally met for an LQ game to be potential. 

\begin{rem} \label{rem:higherdimension} We expect that extending the result to higher dimensions will not, in general, relax the conditions derived for being a potential game. However, extending the proof beyond the scalar case is challenging, if not very tedious, due to reasoning about the set of equalities that arise. Let us consider $t_f = 2$. Observe that for higher dimensions of states, a set of dynamics parameters might exist such that the condition in~(\ref{eq:Lemma7}) is met for this pair of $(t,\tau)$ even though $Q_2^1 \neq Q_2^2$. Such cases might happen when $B^i(B^i)^T$ is not full rank, and consequently, the following condition derived through~(\ref{eq:Lemma7}) for the pair of $(t, \tau) = (1, 0)$ may not have a unique solution.
\begin{flalign*} %\label{eq:HD_condition}
    (B^2)^TA^TQ_{2}^1B^1 = (B^2)^TA^TQ_{2}^2B^1.
\end{flalign*}
In summary, in higher dimensions,  it might be possible to derive slightly relaxed, but system parameter-dependent ($A, B^i$) conditions equivalent to those in Proposition \ref{prop:coupled_i}. 
\end{rem}

\begin{rem}
To investigate the effect of the noisy dynamics on the results of Proposition~\ref{prop:coupled_i}, we add a random noise $w_t$ drawn from the distribution $\mathcal{W}$ with zero mean to the game dynamics in equation (\ref{eq:accum_dynamic}). The loss function will be defined as follows.
\begin{flalign} \label{eq:cost_lqg_noise}
	J^i(\gamma) =& \mathop{\mathbb{E}}_{x_0 \sim \mathcal{D}}\bigl[ \sum_{t = 0}^{t_f-1}\mathop{\mathbb{E}}_{w_t \sim \mathcal{W}} c_t(x_{t+1}, u_t)\bigr],    
\end{flalign}
where the stage loss $c_t(x_{t+1}, u_t)$ is defined the same as the one for the game without noise. Considering the loss function defined above, the definition of Nash equilibrium, potential game, identical interest games for LQ games (Definitions~\ref{def:ne}, \ref{def:potential}, and \ref{def:identical_interest}),  the game \emph{pseudo-gradient}, $\mathcal{G}(\gamma)$, and its Jacobian in equation~(\ref{eq:jac_j}) do not change. Hence, Lemma~\ref{Lm:Lemma7} is still valid for the game with noisy dynamics. Note that the derivative of $c_t(x_{t+1}, u_t)$ with respect to $x_t$ or $u_t^j$ is similar to the game without noise. For open-loop information structure, one can show that the derivatives of $x_{t+1}$ and $u_t^i$ with respect to the decision variable $u_{\tau}^j$ do not depend on $w_t$, which indicates that the results of Proposition~\ref{prop:coupled_i}, part (a) is the same for games with and without noise.  For full-state linear feedback information structure, one can show that the derivatives of $x_{t+1}$ and $u_t^i$ with respect to the decision variable $K_{\tau}^j$ are zero for $t+1 \ge \tau$ and for $t+1<\tau$ can be written as a product of a coefficient that do not depend on $w_t$ and $x_{\tau}$. This observation implies that by investigating conditions in Lemma~\ref{Lm:Lemma7}, the potential conditions derived for the LQ game with noisy dynamics and full-state feedback information structure are the same as those for the game without noise in Proposition~\ref{prop:coupled_i}, part (b). Details are provided in Appendix~\ref{sec:additive_noise}.
\end{rem}

The challenge of defining a class of potential games beyond identical interest games, coupled with the desirable properties of potential games (as discussed in the introduction), motivates us to restrict the game setting. Given the properties of Markov potential games with decoupled dynamics and local policies in finite state and action spaces~\cite{Zhang23}, we confine the LQ game setting, where the state and action spaces are infinite, to decoupled dynamics and the decoupled state linear feedback information structure. Moreover, restricting the $A$ and $B^i$ matrices in LQ games might lead to a potential game that satisfies the conditions in Remark~\ref{rem:higherdimension} without being an identical interest game. As described in Section~\ref{sec:preliminaries}, LQ games with decoupled dynamics are a practically relevant subset of LQ games. On the other hand, the closed-loop Nash equilibrium of the general-sum LQ game, and consequently the LQ game with decoupled dynamics, has a linear full-state feedback structure~\cite{Basar98}. To relax the conditions described in Proposition~\ref{prop:coupled_i}, an information structure sparser than full-state linear feedback, such as decoupled state linear feedback, should be investigated. Thus, we focus on LQ games with decoupled dynamics and a decoupled state linear feedback information structure, similar to the study~\cite{Zhang23} on Markov potential games with decoupled dynamics and localized policies.
% The challenge of defining a class of potential games beyond identical interest on one hand and the desirable properties of a potential game on the other hand (see motivations discussed in the introduction) motivate us to confine the game setting to the decoupled dynamics defined in~(\ref{eq:dyn_decoupled}) and a decoupled state feedback information structure in (iii). In this setting, we identify LQ potential games that are not identical interest games.

\subsection{Decoupled dynamics and  information structure}
\label{sec:potential_decoupled}
Here, we provide sufficient conditions for the LQ game with decoupled dynamics in~(\ref{eq:dyn_decoupled}) and decoupled information structure (iii) to be potential. Furthermore, we show that there exist practically relevant classes of LQ games that satisfy these conditions.   
\begin{thm}\label{Th:CL_po}
	The LQ game defined in~(\ref{eq:cost_lqg}) with decoupled dynamics in~(\ref{eq:dyn_decoupled}) and the decoupled state feedback information structure in (iii) is potential if the following two conditions are satisfied for $i,j \in [N]$ and $t \in \{0, \dots, t_f -1\}$.
	\begin{enumerate}
		\item [(C1)] $(Q^i_{t+1})_{ij} = (Q^j_{t+1})_{ij} : = Q_{t+1}^{ij}$, where $(Q^i_{t+1})_{hl} \in \mathbb{R}^{n^h \times n^l}$ is a sub-matrix of $Q_{t+1}^i$ consisting of the entries relevant to $i$ and $j$ agents as follows
\begin{equation} \label{eq:q_definition}
	Q^i_{t+1} = 
	\begin{bmatrix}
		(Q^i_{t+1})_{11} & \dots & (Q^i_{t+1})_{1N}  \\
                \vdots & \ddots & \vdots \\
		(Q^i_{t+1})_{N1} & \dots&  (Q^i_{t+1})_{NN}
	\end{bmatrix}.
\end{equation} 
		\item [(C2)] $(R_t^i)_{ij} = (R_t^j)_{ij} : = R_t^{ij}$.
	\end{enumerate}
\end{thm}
For a clear comparison of the potential conditions in Proposition~\ref{prop:coupled_i} and Theorem~\ref{Th:CL_po}, let us write the action part of the loss function in equation (\ref{eq:Loss_action_gen}) in a closed form, $J^i_u(\gamma,x_0) =  \sum_{t = 0}^{t_f-1} (u_t)^T R^i_tu_t$ where
\begin{equation}\label{eq:R_definition}
    % \label{eq:potential_FH_parameters}
		 R^i_t :=\begin{bmatrix}
			(R^i_t)_{11} &  \dots & (R^i_t)_{1N}\\
			\vdots &  \ddots & \vdots\\
			(R^i_t)_{N1} &  \dots & (R^i_t)_{NN}
		\end{bmatrix},
\end{equation}
and $u = (u^1, \dots, u^N)\in \mathbb{R}^{N}$. Note that the state part of the loss function $J^i_x(\gamma,x_0)$ in equation (\ref{eq:Loss_state_gen}) can also be written in a closed form when the dynamics are decoupled using the definition of $Q_t^i$ in equation (\ref{eq:q_definition}). Proposition~\ref{prop:coupled_i} states that the game is potential with a full-state feedback information structure if and only if all loss function parameters for both players are the same (except for those corresponding to $t = 0$), i.e., $R^i_t = R^j_t$ and $Q^i_t = Q^j_t$. However, Theorem~\ref{Th:CL_po} indicates that the game is potential with a decoupled state feedback information structure if the non-diagonal elements of $R^i_t$ and $Q^i_t$ are equal to the non-diagonal ones of $R^j_t$ and $Q^j_t$, respectively. It does not impose any conditions on the diagonal elements. As mentioned in Remark~\ref{rem:pastLQpotential}, the LQ games defined in~\cite{Aprem19,Zazo16} with non-zero $(R^i_t)_{ii}$ for $i \in [N]$ are assumed to be potential under full-state linear feedback structure if $Q_t^i = Q_t^j$ and non-diagonal elements of $R^i_t$ and $R^j_t$ are zero.
Although the discussion in Remark~\ref{rem:pastLQpotential} indicates that this claim is wrong, according to Theorem~\ref{Th:CL_po} and assuming decoupled dynamics, one can show that these games are an LQ potential game with decoupled state linear feedback information structure.

% Intuitively, this theorem shows that for two arbitrary agents $i,j$, only the loss parameters corresponding to each other's state and action should be the same. In contrast, in Proposition~\ref{prop:coupled_i} part (b), we required that across all agents, the input and state loss parameters be the same (except for those corresponding to $t = 0$).  

%, only the cost terms corresponding to  Before providing the proof, let us highlight the difference between conditions presented in Proposition~\ref{prop:coupled_i} and Theorem~\ref{Th:CL_po}. By considering conditions (C1) and (C2), the constraints on the loss function parameters of the LQ game are relaxed compared to those in Proposition~\ref{prop:coupled_i}. The first reason is that conditions (C1) and (C2) do not generally satisfy the conditions in Proposition~\ref{prop:coupled_i} part (a) since $Q^i_t$ does not need to be $Q^j_t$. The second reason is that conditions (C1) and (C2) do not generally represent the conditions in Proposition~\ref{prop:coupled_i} part (b) since (C1) and (C2) are not equivalent to the condition in~(\ref{eq:prop5_b}). For example, $(R^i_t)_{lh}$ is not necessarily equal to $(R^j_t)_{lh}$ where $l,h \in [N]/\{i,j\}$.
\begin{pf} The proof leverages Lemma \ref{Lm:Lemma7}. In particular, we verify that an LQ game with decoupled dynamics that satisfies conditions (C1) and (C2) will satisfy the conditions in~(\ref{eq:Lemma7}), and therefore, is potential. Then, by providing more details on conditions, we show that (C1) and (C2) are not necessary for the games to be potential. 

% lead to the biggest subset of potential LQ games with decoupled dynamics and decoupled state feedback information structure
% First, the second derivative of $J_x^i(\gamma,x_0)$ and $J_u^i(\gamma,x_0)$ are derived with respect to .
For the decoupled state feedback information structure in (iii), the decision variable is $\gamma_t^i = k_t^i$ where $u_t^i = -k^i_tx^i_t$. Since the dynamics are decoupled as equation~(\ref{eq:dyn_decoupled}), $u_t^i$ is the function of $(k^i_0, \dots, k^i_t)$; $x_t^i$ is the function of $(k^i_0, \dots, k^i_{t-1})$; and both are independent of the other agents' decision variable, $k^{-i}$. Consequently, for  $i \in [N]$, $J_x^i(\gamma,x_0)$ can be written as follows:
\begin{flalign*} %\label{eq:J_x}
    J_x^i(\gamma,x_0) &=\sum_{t = 0}^{t_f}\sum_{j = 1}^{N}\sum_{h = 1}^{N} \big[ (\Bar{x}^j_{t})^T(Q^i_t)_{jh}\Bar{x}^{h}_{t} \big],
% \\ \label{eq:J_u}
    % J_u^i(k,x_0)&:= \sum_{t = 0}^{t_f-1}\big(\sum_{j = 1}^{N} u_t^jR^{ji}_t\big)u_t^i,
\end{flalign*}
where $\Bar{x}^j_{t} = x^j_t - d^j_t$. For the state part of the loss function, the first derivative of $J^i_x(\gamma,x_0)$ with respect to $k^i_t$ for $t \in \{0, \dots, t_f-1\}$ is computed using the chain rule and product rule as follow,
\begin{flalign*} %\label{eq:Jx_derivative}
        \frac{\partial J_x^i(\gamma,x_0)}{\partial k^i_{t}} = 
        &\sum_{t' = t+1}^{t_f}\sum_{h= 1}^{N} 2(\Bar{x}_{t'}^{h})^T(Q^i_{t'})_{hi}\frac{\partial \Bar{x}_{t'}^i}{ \partial k^i_{t}}.
\end{flalign*}
The second derivatives with respect to $k^i_t$ and $k^j_{\tau}$ are computed by applying the chain rule and product rule such that for $j \in [N] \backslash \{i\}$, and for $\tau \in \{0, \dots, t_f-1\}$,
\begin{flalign*} 
    \frac{\partial^2 J_x^i(\gamma,x_0)}{\partial k_{\tau}^j\partial k_{t}^i} = &\sum_{t' = \text{max}(t,\tau) + 1}^{t_f} 2(\frac{\partial \bar{x}^i_{t'}}{\partial k_{t}^i})^T(Q_{t'}^i)_{ji}^T\frac{\partial \Bar{x}_{t'}^j}{\partial k_{\tau}^j}.
\end{flalign*}

For the action part of the loss function in~(\ref{eq:Loss_action_gen}), similar steps are taken to compute the second derivatives of $J_u^i(\gamma,x_0)$ with respect to $k^i_t$ and $k^j_{\tau}$ as follows 
\begin{flalign*} 
    \frac{\partial^2 J_u^i(\gamma,x_0)}{\partial k_{\tau}^j\partial k_{t}^i} = &  \sum_{t' = \text{max}(t,\tau)}^{t_f-1}  2(R_{t'}^i)_{ji}(\frac{\partial u_{t'}^i}{\partial k_{t}^i})^T\frac{\partial u^j_{t'}}{\partial k_{\tau}^j},
\end{flalign*}
where for $j \in [N] \backslash \{i\}$, and for $\tau \in \{0, \dots, t_f-1\}$. Therefore, the second derivative of $J^i(\gamma)$ with respect to $k^i_t$ and $k^j_{\tau}$ is as follows
\begin{flalign*} 
    \frac{\partial^2 J^i(\gamma)}{\partial k_{\tau}^j\partial k_{t}^i} = &  2\mathop{\mathbb{E}}_{x_0 \sim \mathcal{D}} \Bigg[\sum_{t' = \text{max}(t,\tau)}^{t_f-1}\{(\frac{\partial \bar{x}^i_{t'+1}}{\partial k_{t}^i})^T(Q_{t'+1}^i)_{ji}^T\\
    &\frac{\partial \Bar{x}_{t'+1}^j}{\partial k_{\tau}^j} +(R_{t'}^i)_{ji}(\frac{\partial u_{t'}^i}{\partial k_{t}^i})^T\frac{\partial u^j_{t'}}{\partial k_{\tau}^j}\}\Bigg],
\end{flalign*}
If $(R_{t'}^i)_{ji}$ and $(Q_{t'+1}^i)_{ji}$ are respectively replaced by $(R_{t'}^j)_{ji}$ and $(Q_{t'+1}^j)_{ji}$ in the equation above, the second derivatives of $J^j(\gamma)$ with respect to $k^i_t$ and $k^j_{\tau}$ is achieved. Since the game setting satisfies (C1) and (C2), these loss parameters are equal, i.e. $(R_{t'}^i)_{ji} = (R_{t'}^j)_{ji}$ and $(Q_{t'+1}^i)_{ji}=(Q_{t'+1}^j)_{ji}$, for $t' \in \{0, \dots, t_f-1\}$ and for $i,j \in [N]$. Accordingly, the second derivatives of $J^j(\gamma)$ and $J^i(\gamma)$ are equal, and consequently, the LQ game is potential. 

Next, we show that conditions (C1) and (C2) are not necessary conditions for the game to be potential. For the pair $(t, \tau) = (t_f-1, t_f-1)$, the second derivative of $J^i(\gamma)$ with respect to $k^i_t$ and $k^j_{\tau}$ is
\begin{flalign*} 
    \frac{\partial^2 J^i(\gamma)}{\partial k_{t_f-1}^j\partial k_{t_f-1}^i} = &  2\big\{(b^i)^T(Q_{t_f}^i)_{ij}b^j + (R_{t_f-1}^i)_{ij} \big\}\\
    &\times \mathop{\mathbb{E}}_{x_0 \sim \mathcal{D}} [x_{t_f-1}^i\big(x_{t_f-1}^j\big)^T],
\end{flalign*}
since the derivatives of state and action of agent $i$ are $\frac{\partial \bar{x}^i_{t_f}}{\partial k_{t_f-1}^i} = -b^i(x^i_{t_f-1})^T$ and $\frac{\partial {u}^i_{t_f-1}}{\partial k_{t_f-1}^i} = -(x^i_{t_f-1})^T$. For $(t, \tau) = (t_f-1, t_f-2)$, the second derivative of $J^i(\gamma)$ with respect to $k^i_t$ and $k^j_{\tau}$ is
\begin{flalign*} 
    \frac{\partial^2 J^i(\gamma)}{\partial k_{t_f-2}^j\partial k_{t_f-1}^i} =   2  \big\{(b^i)^T(Q_{t_f}^i)_{ij}A^jb^j - \big[(b^i)^T(Q_{t_f}^i)_{ij}b^j\\
      + (R_{t_f-1}^i)_{ij}\big]k_{t_f-1}^jb^j \big\} \mathop{\mathbb{E}}[x_{t_f-1}^i\big(x_{t_f-1}^j\big)^T],
\end{flalign*}
since the derivatives of states and action of agent $j$ are $\frac{\partial \bar{x}^j_{t_f}}{\partial k_{t_f-2}^j} = -(A^j - b^jk^j_{t_f-1})b^j(x^j_{t_f-2})^T$ and $\frac{\partial {u}^j_{t_f-1}}{\partial k_{t_f-2}^j} = k_{t_f-1}b^j(x^j_{t_f-2})^T$. From the pair of time steps $(t_f-1, t_f-1)$, the condition in~(\ref{eq:Lemma7}) is met if and only if
\begin{flalign*}
    (b^i)^T(Q_{t_f}^i)_{ij}A^jb^j =& (b^i)^T(Q_{t_f}^j)_{ij}A^jb^j, \\
    (b^i)^T(Q_{t_f}^i)_{ij}b^j + (R_{t_f-1}^i)_{ij} =& (b^i)^T(Q_{t_f}^j)_{ij}b^j + (R_{t_f-1}^j)_{ij}.
\end{flalign*}
As discussed in Remark~\ref{rem:higherdimension}, this set of equations may not have a unique solution for $x^i$ with a dimension higher than one. Note that the term $(b^i)^T((A^i)^{t'})^T(Q_{t_f}^i)_{ij}(A^j)^{t"}b^j$ where $t', t" \in \{0, \dots, t_f-2\}$ will appear in several equations if other pairs of $(t, \tau)$ are calculated. This notion indicates that conditions on $A^i$ and $b^i$ might be very restrictive such that the LQ game with decoupled dynamics and decoupled feedback information structure satisfying conditions (C1) and (C2) might be the only game that meets condition~(\ref{eq:Lemma7}) for all pairs of $(t, \tau)$.\qed
\end{pf}
% Thus, the LQ game in Theorem~(\ref{Th:CL_po}) is a big subset of potential LQ games with decoupled dynamics and decoupled feedback information structure.

\begin{rem}
To investigate the effect of the noisy dynamics on the results of Theorem~\ref{Th:CL_po}, we add a random noise $w^i_t$ drawn from the distribution $\mathcal{W}^i$ with zero mean to the set of decoupled dynamics in equation (\ref{eq:dyn_decoupled}). By defining $w_t := (w^1_t, \dots, w^N_t)$, the loss function of the LQ game with noisy decoupled dynamics can be written as the one in equation~\ref{eq:cost_lqg_noise}. Thus, Lemma \ref{Lm:Lemma7} is still valid for the game with noisy dynamics and decoupled state linear feedback information structure, and the derivative of $c_t(x_{t+1}, u_t)$ with respect to $x_t$ or $u_t^j$ is similar to the game without noise. For decoupled state linear feedback information structure, the derivatives of $x^i_{t+1}$ and $u_t^i$ with respect to the decision variable $k_{\tau}^i$ are zero for $t+1 \ge \tau$ and for $t+1<\tau$ can be written as a product of a coefficient that do not depend on $w_t$ and $x^i_{\tau}$. Note that the derivatives of $x^i_{t+1}$ and $u_t^i$ with respect to $k_t^j$ for $j \neq i$ are zero. This observation implies that by investigating conditions in Lemma \ref{Lm:Lemma7}, the potential conditions derived for the LQ game with noisy decoupled dynamics and full-state feedback information structure are the same as those for the game without noise in Theorem~\ref{Th:CL_po}.  Details are provided in Appendix~\ref{sec:additive_noise}.
\end{rem}
After theoretically deriving non-trivial (non-identical interest) classes of state-feedback potential games, we confirm the existence of practically relevant problems that meet the conditions in Theorem~\ref{Th:CL_po}.

\textbf{Example 1 continued:}  In the formation control games in Example 1, condition (C2) is satisfied since $(R^i_t)_{ij} = 0$ for $i \in [N]$ and $j \in [N] \backslash \{i\}$ . If the weights are symmetric such that $w^{ij}_t = w^{ji}_t$, the state loss parameter $(Q^i_t)_{ij} = w^{ij}_t$ will be equal to $(Q^j_t)_{ij}=w^{ji}_t$; hence, condition (C1) is also met. Thus, from Theorem~\ref{Th:CL_po}, Example 1 with decoupled feedback information structure is an LQ potential game.

\textbf{Example 2 continued:} The dynamic Cournot game described in Example 2 with the loss function in~(\ref{eq:cournot_loss}) and the decoupled feedback information structure is also potential. The reason is that this game meets conditions (C1) and (C2) in Theorem~\ref{Th:CL_po} as $(R^i_t)_{ji} = p_{t}$, $Q^i_t = 0$ for $t<t_f$, and $(Q^i_{t_f})_{ij} = 0$ where $i \in [N]$ and $j \in [N] \backslash \{i\}$.

% represents the global objective in the formation control games. In other words, the terms in the loss function of an agent representing its distance from another agent should be similar in the other agent loss function.

% with decoupled dynamics and information structure (iii) is shown to be potential. This game is not an identical interest game, as $(Q^i_{t_f})_{ii}$ need not be $(Q^j_{t_f})_{ii}$.
% \begin{cor}
%     The dynamic Cournot game described in example (2) with the loss function defined in~(\ref{eq:cournot_loss})
%     % \begin{align*}
%     %     J^i(x_0,u) = \sum_{t=0}^{t_f -1} \left( p_{1,t} \sum_{j=1}^N u^j_t - p_{0,t} \right) u^i_t + Q^i_{tf} ( x^i_{t_f} - \overline{x}^i)^2,
%     % \end{align*}
%     and the decoupled state feedback information structure in (iii) is potential.
% \end{cor}
% \textit{Proof}:
% 	The dynamic Cournot game in Example (2) has a decoupled dynamic, and with the decoupled state feedback information structure in (iii), it meets conditions (C1) and (C2) in Theorem~(\ref{Th:CL_po}) as $(R^i_t)_{ji} = p_{1,t}$, $Q^i_t = 0$ for $t<t_f$, and $(Q^i_{t_f})_{ij} = 0$ where $i \in [N]$ and $j \in [N] \backslash \{i\}$. \qed

What is the potential function of a game in the class identified above? Does a Nash equilibrium exist for such a game? If yes, how could we compute it? These are the questions we examine in the next section. While answers are provided to the first two questions, we identify the challenges in answering the third one. 

\section{Properties of decoupled LQ potential games}\label{decoupled_game}
We first derive the potential function for the game satisfying the conditions of Theorem~\ref{Th:CL_po}. This potential function corresponds to a single-agent optimal control problem. However, this problem is not a standard LQ control due to the structure of the control policy arising from the decoupled state linear feedback information structure. Given that any optimizer of the potential function will be a Nash equilibrium, we verify that the control problem has an optimal solution.  However,  we show that characterizing or computing the best response is non-trivial. 

%Then, we prove that there is at least one Nash equilibrium.  show that the best response for each agent is equivalent to the local loss's best response if the agent's initial states are not correlated. Next, we explain that the stationary point might not be a Nash equilibrium.

\subsection{Existence of a Nash equilibrium}
\begin{prop} \label{prop:pot_function}
    The LQ game in~(\ref{eq:cost_lqg}) with decoupled dynamics in~(\ref{eq:dyn_decoupled}) and the decoupled state feedback information structure in (iii) under conditions (C1) and (C2) in Theorem~\ref{Th:CL_po} has the potential function $\Pi: \mathbb{R}^{t_fn} \rightarrow \mathbb{R}$: 
    \begin{equation} \label{eq:LQR_loss}
	   \Pi (\gamma)=\mathop{\mathbb{E}}_{x_0 \sim \mathcal{D}} \big[\sum_{t =0}^{t_f} \bar{x}_{t}^TQ_{t}\bar{x}_{t} +
       \sum_{t=0}^{t_f-1} u_t^TR_tu_t\big],
	\end{equation}
    where $\bar{x}_{t} = x_t - d_t$, and the joint state, control, and decision variable are respectively, $x_t$, $u_t$, and $\gamma = k$. The loss function matrices $R_t \in \mathbb{R}^{N \times N}$, and $Q_t \in \mathbb{R}^{n \times n}$ are
    \begin{equation*}
    % \label{eq:potential_FH_parameters}
		 R_t =\begin{bmatrix}
			R_t^{11} &  \dots & R_t^{1N}\\
			\vdots &  \ddots & \vdots\\
			R_t^{N1} &  \dots & R_t^{NN}
		\end{bmatrix}, \ Q_t =\begin{bmatrix}
			Q_t^{11} &  \dots & Q_t^{1N}\\
			\vdots &  \ddots & \vdots\\
			Q_t^{N1} &  \dots & Q_t^{NN}
		\end{bmatrix},
	\end{equation*}
    with $Q^{ij}_t$, $R_t^{ij}$ defined in (C1) and (C2) in Theorem~\ref{Th:CL_po}. Note that $R_t$ is positive definite, and $Q_t$ is positive semi-definite.
\end{prop} 
\begin{pf}
The LQ game with decoupled dynamics and the decoupled state feedback information structure in (iii) is a potential game under conditions (C1) and (C2) in Theorem~\ref{Th:CL_po}. For this LQ game, we illustrate that the loss function of agent $i$ is a summation of equation~(\ref{eq:LQR_loss}) and a dummy game which is defined in~\cite{Hespanha17} as a game that is independent of agent $i$ decision variables $k^i$, and according to Proposition 12.5 in~\cite{Hespanha17}, the function in equation~(\ref{eq:LQR_loss}) is the potential function of this LQ game.

For the state part, in the definition of $Q_t$, if $Q_t^{hl}$ where $h,l \in [N] \backslash \{i\}$ is replaced by $(Q_t^i)_{hl}$, the matrix $Q^i_t$ is formed, and consequently, the state part of the $i$th agent loss function at each time step, $(x_t-d_t)^TQ_t^i(x_t-d_t)$, is achieved. Note that $(x^h_t-d^h_t)^TQ_t^{hl}(x^h_t-d^h_t)$ and $(x^h_t-d^h_t)^T(Q^i_t)_{hl}(x^h_t-d^h_t)$ are independent of the agent $i$'s decision.  This means replaced terms were dummy games for agent $i$.
For the action part, in the definition of $R_t$, if $R_t^{hl}$ where $h,l \in [N] \backslash \{i\}$ is replaced by $(R_t^i)_{hl}$, $u_t^TR_tu_t$ becomes equal to the action part of $i$th agent loss function at time $t$. Note that $u_t^hR^{hl}_tu_t^l$ and $u_t^h(R_t^i)_{hl}u_t^l$ where $h,l \in [N] \backslash \{i\}$ are independent of agent $i$'s policy and represent dummy games. \qed 
\end{pf}
%The details of the proof are provided in Appendix \ref{ap:pot_function}, and here, the intuition behind the proof is stated. 
% The key steps in the proof involve forming a new loss function matrices $Q_t$ and $R_t$ by bringing the entries from agents' loss functions, which correspond to the quadratic terms and the cross-terms that are similar among all agents.
In contrast to a classical control problem, the controller must maintain a certain sparsity due to the information structure defined in (iii). Thus, we need to examine the existence of an optimal solution, which corresponds to a Nash equilibrium of the game. 
%Moreover, one can verify that it is not quadratic invariant, thus the results obtained in \cite{Furieri} can not be applied. Thus we necessitating a foundational examination of the problem. 
In Proposition \ref{prop:NEminimum} below, we address the existence of Nash equilibrium and discuss the challenge in its characterization in Proposition \ref{prop:initialStatesCorrelation}. 
%The first question we should answer is whether or not a Nash equilibrium exists. In the case of a potential function, this is equivalent to determining if a global optimum exists. We now prove the existence of at least one Nash equilibrium, and later, we discuss the open questions.

\begin{prop} \label{prop:NEminimum}
    Any LQ potential game with decoupled dynamics and the decoupled state feedback information structure in (iii) has at least one Nash equilibrium. This  Nash equilibrium is the optimal structured control of the single-agent loss function identified in~(\ref{eq:LQR_loss}).
\end{prop}
\begin{pf}
 The loss function is continuous and bounded below, growing to infinity whereas any $K_t$ grows to infinity. Thus, at least an optimal control policy exists. This policy is a Nash equilibrium as any potential function optimizer is a Nash equilibrium of the potential game. \qed
\end{pf}

Proposition \ref{prop:NEminimum} does not provide insight into the computation of a Nash equilibrium.
To address this, we first show that the agents' best responses are coupled if and only if their initial states are correlated (see Proposition \ref{prop:initialStatesCorrelation}). In other words, with uncorrelated initial states, the game becomes degenerate, allowing agents to compute the optimal policies using only their local loss function and dynamics. Thus, the only non-trivial LQ potential games with decoupled dynamics and losses are the ones with correlated initial states. However, we also highlight the challenge in computing a Nash equilibrium in this case by demonstrating that dynamic programming cannot be used to compute the best response, and therefore, a Nash equilibrium of the game.  
%Before delving into this topic, we must first introduce two propositions. These propositions demonstrate that, when the initial states are correlated, the optimal policies of agents become coupled, and dynamic programming cannot be employed to compute an agent's best response.
 \begin{prop} \label{prop:initialStatesCorrelation}
 The best response of an agent $i$ in the LQ game defined in~(\ref{eq:cost_lqg}) with decoupled dynamics in~(\ref{eq:dyn_decoupled}) and the decoupled state feedback information structure in (iii) under conditions (C1) and (C2) in Theorem \ref{Th:CL_po} depends on agent $j$'s state parameters and initial state if and only if their initial states are correlated.
 \end{prop}
\begin{pf}
The loss function is bounded below; thus, at least one best response policy exists. If there is a unique value for which the derivative equals zero, the best response can be determined by finding the stationary point. By taking the derivative of the loss function of agent $i$ with respect to $k^i_\tau$, we obtain:
\begin{align*}
    &\frac{\partial J^i(k^i, k^{-i})}{\partial k^i_\tau} = 2 \sum_{t = \tau + 1}^{t_f} \sum_{j = 1}^{N} (b^i)^\top \bigl( (Q_t^i)_{ij} + (k^i_t)^\top (R^i_t)_{ij} k^j_t \bigr) \\
    &\mathop{\mathbb{E}}_{x_0 \sim \mathcal{D}} \left[ x^j_t (x^i_{t,-\tau})^\top \right] + 2 \sum_{j = 1}^{N} (R^i_\tau)_{ij} k^j_\tau \mathop{\mathbb{E}}_{x_0 \sim \mathcal{D}} \left[ x^j_\tau (x^i_\tau)^\top \right],
\end{align*}
where $x^i_{t,-\tau} = \prod_{t'= 0, t' \neq \tau}^t (A^i - b^i k^i_{t'} ) x_0^i$. 

First, we show that the correlation of the initial states of two agents is a necessary condition for their policies to depend on each other's parameters and initial states. If the initial state of agent $i$ is not correlated with any other agent's initial state, it indicates that $ \mathop{\mathbb{E}} \left[ x_0^j (x_0^i)^\top \right] = 0$ for $j \neq i$, and  $ \mathop{\mathbb{E}} \left[ x_\tau^j (x_\tau^i)^\top \right] = 0$ for $j \neq i$ and any $\tau \in \{0,\dots, t_f\}$, because $\mathop{\mathbb{E}} \left[ x_\tau^j (x_\tau^i)^\top \right] = \prod_{t' = 1}^{\tau - 1} (A^j - b^j k^j_{t'} ) \mathop{\mathbb{E}} \left[ x_0^j (x_0^i)^\top \right] \left( \prod_{t' = 1}^{\tau - 1} (A^i - b^i k^i_{t'} ) \right)^\top $. Hence, the best response $k^{i,*}_{t_f-1}$ is:
\begin{align*}
    k^{i,*}_{t_f-1} =  \left( (R^i_{t_f-1})_{ii} + (b^i)^\top (Q_{t_f}^i)_{ii} b^i \right)^{-1} (b^i)^\top (Q_{t_f}^i)_{ii} A^i.
\end{align*}
The other $k^{i,*}_\tau$, for $\tau \in \{0, \dots, t_f - 2\}$, can be computed through dynamic programming and they do not depend on other agents' parameters as well. Thus, the best response for agent $i$ depends only on its own parameters and initial state.

Now, we show that the correlation of the initial states of two agents is a sufficient condition for their policies to depend on each other's parameters and initial states. For the sake of simplicity, we take $t_f = 1$, but the same results will hold for any value of $t_f$. Assuming that $\mathop{\mathbb{E}} \left[ x_0^i (x_0^i)^\top \right]$ is invertible, there exists a unique stationary point and a unique $k^{i,*}_0$:
\begin{align*}
    k^{i,*}_0 =& \left( (R^i_{0})_{ii} + (b^i)^\top (Q_{1}^i)_{ii} b^i \right)^{-1} \biggl\{ (b^i)^\top (Q_{1}^i)_{ii} A^i + \\
    & \sum_{j \neq i} \bigl( (b^i)^\top (Q_1^i)_{ij} (A^j - b^j k^j_0 ) - \\ 
    & (R^i_0)_{ij} k^j_0 \bigr) \mathop{\mathbb{E}} \left[ x_0^j (x_0^i)^\top \right] \biggr( \mathop{\mathbb{E}} \left[ x_0^i (x_0^i)^\top \right] \biggl)^{-1} \biggr\}.
\end{align*}
If the correlation between the initial states of agent $i$ and at least one other agent $j$ is non-zero, the optimal $k^{i,*}_{0}$ depends on the value of $A^j$, $b^j$, $k^j$, and $x^j_0$. This result holds even for non-potential, decoupled LQ games. 
\qed
\end{pf}

The proof above can be generalized to show that the computation of the optimal decoupled linear control policy of an agent $i$ using dynamic programming is not possible if its initial state is correlated with the initial state of at least one other agent. 
In particular, the best response $k^{i,*}_{t_f-1}$ depends on the values of $x_{t_f-1}^i$, which in turn depends on the policies at the previous time steps. This interdependence precludes the use of the dynamic programming method to compute the best response.
\begin{comment}
\begin{prop} \label{prop:DP}
Under the assumption that the initial state of agent $i$ is correlated with the initial state of at least one other agent, the dynamic programming method cannot compute the best response among decoupled linear policies of the LQ potential game introduced in the Theorem \ref{Th:CL_po}.
\end{prop}
%\begin{pf}
\end{comment}

%\qed    
%\end{pf}
    %It is thus equal to:
%\begin{align*}
    %k^{i,*}_{t_f-1} =&\left( (R^i_{t_f-1})_{ii} + (b^i)^\top (Q_{t_f}^i)_{ii} b^i \right)^{-1} \bigg\{ (b^i)^\top (Q_{t_f}^i)_{ii} A^i  \\
    %&+ \sum_{j = 1, j \neq i}^{N} \biggl( \bigl( (b^i)^\top (Q_{t_f}^i)_{ij} (A^j + b^j k^j_{t_f-1})\\
    %&+ (R^i_{t_f-1})_{ij} k_{t_f-1}^j \biggr) \mathbb{E} \left[ x_{t_f-1}^j (x_{t_f-1}^i)^\top \bigr) \right] \\
    %&\left( \mathbb{E} \left[ x_{t_f-1}^i (x_{t_f-1}^i)^\top \right] \right)^{-1} \bigg\}. 
%\end{align*}
% \qed

%Since the dynamic programming method cannot be employed, a revisitation of the results previously established for the linear quadratic regulator and general-sum LQ games becomes necessary. 

\subsection{Challenge in computing the Nash equilibria}
\label{sec:challenge}
%The results obtained in \cite{Basar98} for the uniqueness of general-sum LQ games Nash equilibrium were based on the dynamic programming induction and, therefore, do not directly apply to our setup. Moreover, due to the complex relationships between agents' optimal control policies and the fact that we could not completely characterize them, obtaining guarantees on the uniqueness (or lack thereof) of the Nash equilibrium proves difficult.
Recent works have proven that policy gradient can be applied for computing the optimal control policy in LQ problems~\cite{Fazel18} and LQ problems with sparsity constraints~\cite{Furieri}. For the LQ regulator problem,~\cite{Fazel18} proved that policy gradient converges to the unique optimum, based on the uniqueness of the global optimum and the gradient-dominated and almost smooth nature of the loss function (see Lemmas 3 and 6 in \cite{Fazel18}). In \cite{Furieri}, this result was extended to finite-horizon structured control problems, a class that encompasses our decoupled LQ potential games. However, the convergence proof in~\cite{Furieri} relied on the so-called quadratic invariance property of the sparsity constraints. Unfortunately, as shown below, the quadratic invariance is not generally satisfied for the decoupled information structure considered here.

\textbf{Counterexample for the quadratic invariance:}
\label{sec:qi}
A subspace $\mathcal{K}$ is quadratic invariant with respect to a matrix $\textbf{CP}_{12}$ if and only if
\begin{align*}
    \textbf{KCP}_{12} \textbf{K} \in \mathcal{K} \text{,  } \qquad \forall \textbf{K} \in \mathcal{K},
\end{align*}
where $\textbf{C}$ and $\textbf{P}$ are matrices depending on the system parameters.

We provide an example of a decoupled LQ game with two agents and time horizon $t_f = 2$ that does not respect the condition of quadratic invariance as presented in~\cite[Definition 1]{Furieri}. The notation used here is aligned with that of~\cite{Furieri}. For simplicity, we choose the system matrices $A$ and $B$ as the identity matrix $A = B = I_2.$
We also assumed perfect knowledge of the state, thus the matrix $C$ is the identity matrix as well.

The subspace $\mathcal{K}$ must satisfy the sparsity constraint due to the information structure outlined in (iii). In our scenario, $\mathcal{K} = \text{diag}(k_1,k_2,k_3,k_4)$ with  $k_1,k_2,k_3,k_4 \in \mathbb{R}$. Given the matrix $\textbf{K} = I_4 \in \mathcal{K}$, we have that $
    \textbf{KCP}_{12} \textbf{K} = \begin{bmatrix}
        &0_{2 \times 2} &0_{2 \times 2}\\
        &I_2 &0_{2 \times 2}
    \end{bmatrix} \notin \mathcal{K}.
$
Thus, the subspace $\mathcal{K}$ is not quadratic invariant with respect to matrix $\textbf{CP}_{12}$. Note that the issue here is that decoupled policies are memoryless. Namely, if we allow the policy of agent $i$ at time $t$ to depend on its state at previous times the same information structure can result in quadratic invariance.

Consequently, even if we can prove the convergence of the policy gradient to a stationary point, this point is not guaranteed to be a Nash equilibrium. Due to the nonlinear dependence of $k^{i,*}_t$ on the policies at previous steps, making assumptions about the uniqueness of the set $(k^{i,*}_0, \dots, k^{i,*}_{t_f-1})$ for which the gradient is zero becomes difficult. The roots of a set of $t_f$ non-linear and non-convex equations may represent multiple stationary points, not all of which are the optimal structured control. Thus, while a Nash equilibrium is a stationary point, other stationary points may exist that are not Nash equilibria.

%The cost function is non-convex, and the proof for the gradient domination property presented in~\cite{Fazel18} does not hold for this game setting. Thus, unlike the results obtained for the linear quadratic regulator in~\cite{Fazel18}, having a null gradient  is not a sufficient condition to assert that the policy is the global optimum.  For these reasons, while we have established the existence of at least one Nash equilibrium, determining specific equilibria poses significant challenges.

\subsection{Convergence of policy gradient  to a stationary point}
Here, we derive conditions for convergence to a stationary point of the policy gradient. Consider the joint decision-variables $k = (k^1, k^2,..., k^N) \in \mathbb{R}^{t_f nN}$, where $k^i = (k^i_1,...,k^i_{t_f-1})$. Starting from an arbitrary $k^i(0) \in \mathbb{R}^{t_fn} $, each agent updates its feedback policy as
\begin{align} \label{eq:updateAgenti}
    k^i(m+1) = k^i(m) - \eta(m) \lambda^i \nabla_{k^i} J^i(k(m)).
\end{align}
The step size of each agent is composed of two terms. The first one, $\eta(m) \in \mathbb{R}_{> 0}$, is the same for all agents, while the second term, $\lambda^i \in \mathbb{R}_{> 0}$, can be chosen independently by each agent. From the definition of derivative and equation (\ref{eq:PotentialDef}), it follows that the gradient of the loss function $J^i(k)$ with respect to $k^i$ is equal to the gradient of the potential function $\Pi (k)$ defined in (\ref{eq:LQR_loss}) with respect to $k^i$. The update iteration~(\ref{eq:updateAgenti}) is equivalent to:
    \begin{equation} \label{eq:gradk_i}
        k^i(m+1) = k^i(m) - \eta(m) \lambda^i \nabla_{k^i} \Pi (k(m)).
    \end{equation}
    Stacking the iterations in~(\ref{eq:gradk_i}) for all agents we obtain:
    \begin{align} \label{eq:gradient}
        k(m+1) &= k(m) - \eta(m) \Lambda \nabla \Pi (k(m)),\\
        \Lambda &= \text{blockdiag}( \lambda_1 I_{t_f n^1}, ..., \lambda_N I_{t_f n^N} \nonumber).
    \end{align}
%    Note that $\Lambda$ is positive definite because all $\lambda_i$ are positive real numbers. 
%    Thanks to the Proposition \ref{prop:NEminimum} we know there is at least one stationary point. In the next theorem, we prove that the policy gradient algorithm proposed converges to a stationary point of the game under a suitable assumption on the step size $\eta(m)$.
\begin{comment}
\begin{assum}
The step size $\eta(m)$ respects the Wolfe conditions for any $k(m)$ and $\Lambda \nabla \Pi(k(m))$.
\end{assum}
\end{comment}
Note that the term $-\Lambda \nabla \Pi (.) \in \mathbb{R}^{t_fn}$ is a descent direction, $
 -(\nabla \Pi(.))^T \Lambda \nabla \Pi(.) < 0,$ and $\Pi$ is continuously differentiable and bounded below in every direction. Thus, from Lemma 3.1 in \cite{Nocedal}, there exists an $\eta$ satisfying the Wolfe conditions, and it can be found using the bisection algorithm presented in~\cite[Proposition 5.5]{Aragon}. It follows from \cite{Nocedal} that   $
        \sum_{m \geq 0} \biggl( \frac{\nabla \Pi(k(m))^T \Lambda \nabla \Pi(k(m))}{ \| \Lambda \nabla \Pi(k(m))\|} \biggr)^2  < \infty.
$
This condition is true if and only if $\nabla \Pi (k_m)$ converges to $0$, which implies that:
% \begin{equation*}
      $  \lim_{m \to \infty} k(m) = k^s,$
    % \end{equation*}
where $k^s$ is a stationary point. Note that variants of the policy gradient above, such as stochastic policy gradient assuming no access to the model, can be developed.

\subsection{Simulation}
\label{sec:simulation}
In this subsection, we implement the policy gradient method introduced in equation~(\ref{eq:gradient}). Our goal is to evaluate the convergence of the policy gradient to a stationary point.

We consider a game with $N = 4$ agents, where each agent has a state of dimension $n^i = 1$ and scalar control actions. The time horizon is $t_f = 10$. To construct the example, we randomly generate the parameters of the game. The values of $ A^i\in \mathbb{R}$ and $b^i \in \mathbb{R}$ are independently sampled from a normal distribution with zero mean and unit variance. Since the matrices $Q_t \in \mathbb{R}^{4 \times 4}$ and $R_t \in \mathbb{R}^{4 \times 4}$ of the potential function should be respectively positive semidefinite and positive definite, we first generate the auxiliary matrices $M_t \in \mathbb{R}^{4 \times 4}$ and $N_t \in \mathbb{R}^{4 \times 4}$ with entries drawn from the same normal distribution, and then we set $Q_t = M_t^\top M_t$ and $R_t = N_t^\top N_t + 0.01 I_4$.

\begin{figure}[tb] 
    \centering
    \begin{subfigure}{\linewidth}
        \centering
        \includegraphics[width=\linewidth]{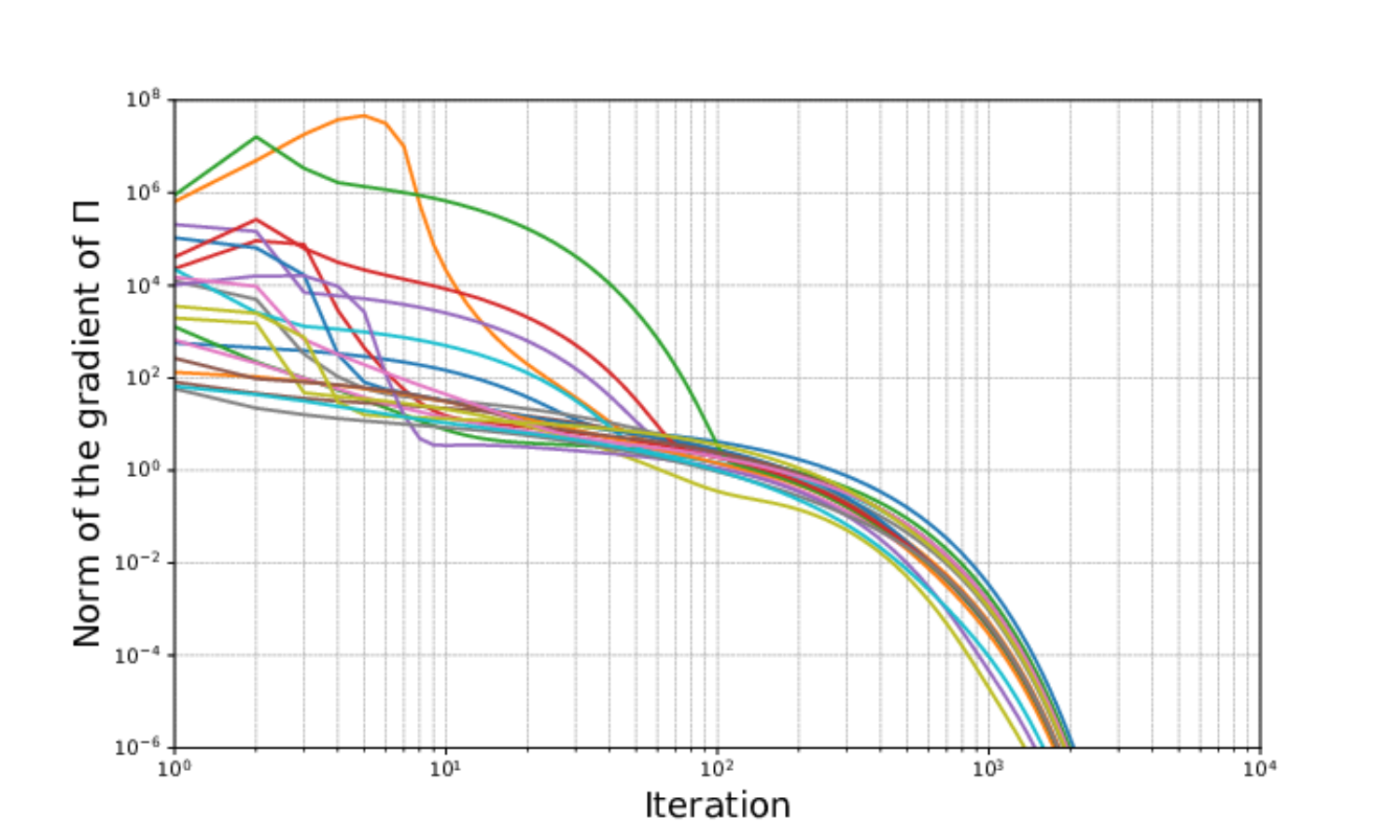}
        \caption{$|\nabla \Pi(k(m))|$}
        \label{fig:grad}
    \end{subfigure}
    \vfill
    \begin{subfigure}{\linewidth}
        \centering
        \includegraphics[width=\linewidth]{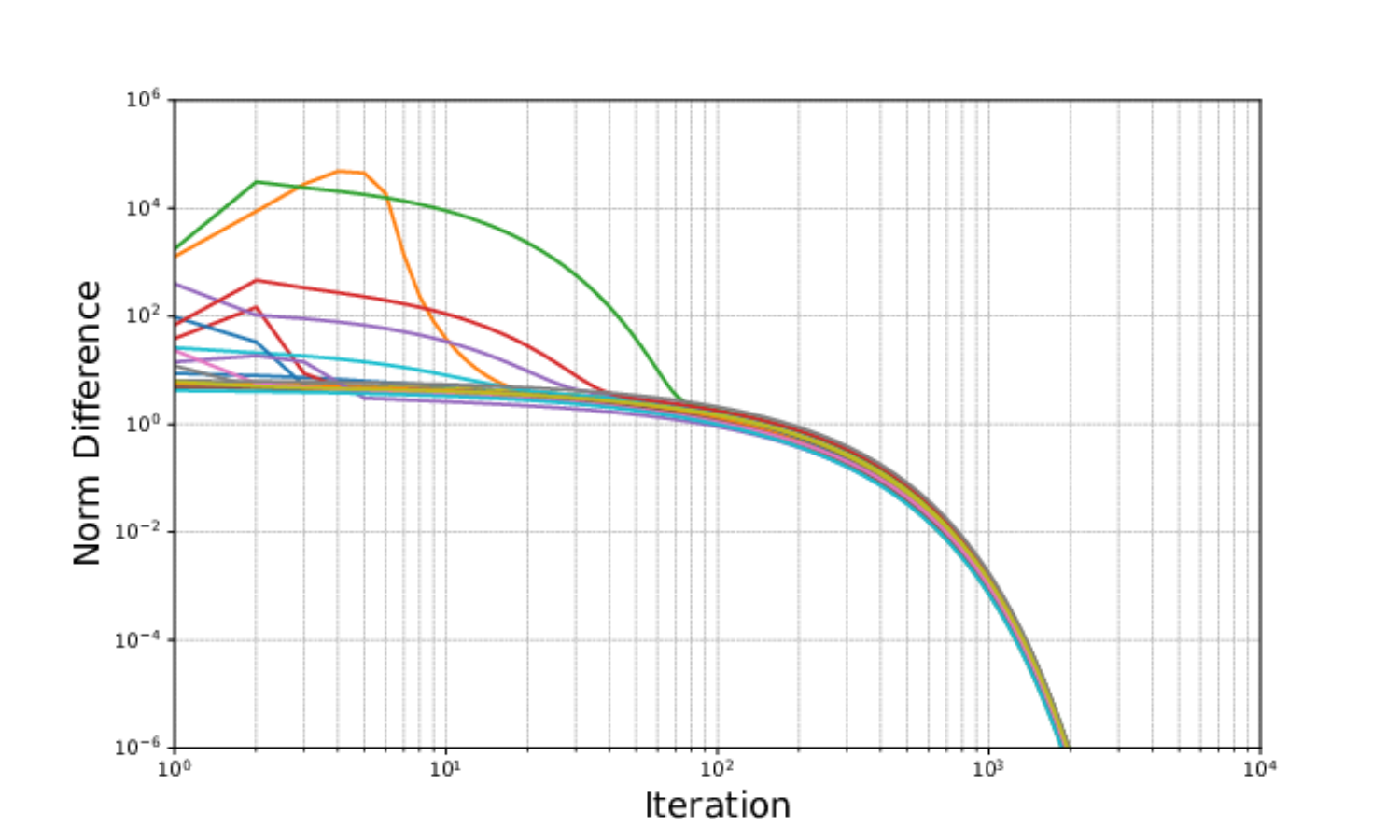}
        \caption{ $|k(m) - \tilde{k}|$}
        \label{fig:K}
    \end{subfigure}
    \caption{Convergence of policy gradient with the update rule in equation~(\ref{eq:gradient}) to a stationary point $\tilde{k}$ for 20 simulations. In these figures, $m$ indicates the iteration number, and the abscissa and the ordinate are on a logarithmic scale.}
\end{figure}

Starting from a random $k(0)$, we run the policy iteration for $10^4$ steps. We observe that $k(m)$ converges to a stationary point $\tilde{k}$. Then, we repeat the experiment twenty times with different random initializations $k(0)$. The values $k^i_t(0)$ are randomly sampled from a standard normal distribution. As shown in Fig~\ref{fig:grad}, all simulations converge to a stationary point. Furthermore, Fig~\ref{fig:K} shows that the distance from the previously computed \( \tilde{k} \) approaches zero in all simulations, indicating convergence to the same stationary point. As proved in Proposition \ref{prop:NEminimum}, the game admits at least a Nash equilibrium. Since we found only one stationary point, and any Nash equilibrium has to be a stationary point, it would seem plausible that the $\tilde{k}$ is a Nash equilibrium. However, as discussed in Subsection \ref{sec:challenge}, proving that this point is a Nash equilibrium requires demonstrating that no agent can improve their outcome by unilaterally selecting a different policy. Due to the nonlinear dependence of the optimal strategy \(k^{i,*}_t\) on the policies at both previous and successive steps, checking for potential improvements would necessitate an exhaustive search over the entire policy space. Specifically, this involves exploring a space with dimensions \(n \times t_f = 40\) in our current setting.

\begin{comment}    
\begin{thm} \label{th:statio}
    For any $k^i(0)$ and $\lambda^i > 0$, $i \in [N]$, the iteration~(\ref{eq:gradient}) will converge to a stationary point.
\end{thm}
\begin{pf}
The proof is based on Theorem 3.2 from \cite{Nocedal}. The potential function $\Pi$ is polynomial in $k$, so the gradient $\nabla \Pi$ will be polynomial. For this reason, it is always Lipschitz continuous in any compact set $D$. Thus, from Theorem 3.2 in \cite{Nocedal} we obtain:
    \begin{equation} \label{eq:convergenceInfty}
        \sum_{m \geq 0} \biggl( \frac{\nabla \Pi(k(m))^T \Lambda \nabla \Pi(k(m))}{ \| \Lambda \nabla \Pi(k(m))\|} \biggr)^2  < \infty.
    \end{equation}
    Finally, the equation~(\ref{eq:convergenceInfty}) is true if and only if $\nabla \Pi (k_m)$ converges to $0$, which implies that:
    % \begin{equation*}
          $  \lim_{m \to \infty} k(m) = k^s,$
        % \end{equation*}
    where $k^s$ is a stationary point.
    \qed
\end{pf}
\end{comment}

%%%%%%%%%%%%%%%%%%%%%%%%%%%%%%%%%%%%%%%%%%%%%%%%%%%%%%%%%%%%%%%%%%%%%%%%%%%%%%%%

\section{Conclusion} \label{sec:Conclusion}
Motivated by the desirable properties of dynamic potential games, our work focused on advancing a fundamental understanding of linear quadratic potential games. To this end, we considered finite-horizon linear quadratic games and derived conditions under which the game would be potential. Notably, we showed that an example of a general class of scalar linear quadratic games with full-state feedback is potential if it deviates slightly from an identical interest game. Additionally, we investigated a subclass of games with decoupled dynamics and decoupled state feedback, revealing the existence of potential games beyond identical interest games. Furthermore, we analyzed the potential function of this subclass and demonstrated the existence of at least one Nash equilibrium. We illustrated how initial state correlations can influence agents' behavior and highlighted the challenge of computing a Nash equilibrium policy. 

Future research can explore approaches to compute the Nash equilibria of the identified potential LQ games, investigate alternative control structures beyond the decoupled state linear feedback, and extend the potential LQ game characterization to an infinite-horizon setting. 
% \textcolor{red}{We also encourage further analysis on more generalized settings of LQ games, for example, including the products of state and action (cross terms) in the loss function, which might lead to more interesting characterization by relaxing the potential conditions derived in Proposition~\ref{prop:coupled_i} and Theorem~\ref{Th:CL_po}.}

\begin{ack}                               % Place acknowledgements
Sara Hosseinirad is financially supported by the Natural Sciences and Engineering Research Council of Canada.  Giulio Salizzoni is financially supported by the Swiss National Science Foundation, grant number 207984. % here.
\end{ack}

\bibliographystyle{plain}        % Include this if you use BibTeX 

% \begin{thebibliography}{10}
% \bibliography{autosam}           % and a bib file to produce the 
                                 % bibliography (preferred). The
                                 % correct style is generated by
                                 % Elsevier at the time of printing.
% \begin{thebibliography}{99}     % Otherwise use the  
                                 % thebibliography environment.
                                 % Insert the full references here.
                                 % See a recent issue of Automatica 
                                 % for the style.
% \input{output.bbl}
%  \bibitem[Heritage, 1992]{Heritage:92}
%     (1992) {\it The American Heritage. 
%     Dictionary of the American Language.}
%     Houghton Mifflin Company.
%  \bibitem[Able, 1956]{Abl:56}
%     B.~C.~Able (1956). Nucleic acid content of macroscope. 
%     {\it Nature 2}, 7--9. 
%  \bibitem[Able {\em et al.}, 1954]{AbTaRu:54}   
%     B.~C. Able, R.~A. Tagg, and M.~Rush (1954).
%     Enzyme-catalyzed cellular transanimations.
%     In A.~F.~Round, editor, 
%     {\it Advances in Enzymology Vol. 2} (125--247). 
%     New York, Academic Press.
%  \bibitem[R.~Keohane, 1958]{Keo:58}
%     R.~Keohane (1958).
%     {\it Power and Interdependence: 
%     World Politics in Transition.}
%     Boston, Little, Brown \& Co.
%  \bibitem[Powers, 1985]{Pow:85}
%     T.~Powers (1985).
%     Is there a way out?
%     {\it Harpers, June 1985}, 35--47.

% \end{thebibliography}

\appendix
\section{Additive Noise} \label{sec:additive_noise}
Here, we prove that the results in Proposition~\ref{prop:coupled_i} and Theorem~\ref{Th:CL_po} can be generalized to the LQ game with noise. For this purpose, we introduce the loss function of the game with noisy dynamics, compare the derivatives in condition~\eqref{eq:Lemma7} of Lemma~\ref{Lm:Lemma7} for this loss function and the one without noise, and illustrate that the final potential conditions in Proposition~\ref{prop:coupled_i}  and Theorem~\ref{Th:CL_po} are the same for both loss functions.
We assume the dynamics of the game include noise as follows.
\begin{equation}
\label{eq:accum_dynamic}
	x_{t+1} = Ax_{t} + \sum_{i = 1}^{N}B^iu^i_{t} + w_t, \quad x_0 \sim \mathcal{D}, \ w_t \sim \mathcal{W}
\end{equation}
where $w_t$ is a random noise that is assumed to be drawn from the distribution $\mathcal{W}$ with zero mean. The loss function will be as follows.
\begin{flalign} \label{eq:cost_lqg_report}
	J^i(\gamma) =& \mathop{\mathbb{E}}_{x_0 \sim \mathcal{D}}\bigl[ \sum_{t = 0}^{t_f-1}\mathop{\mathbb{E}}_{w_t \sim \mathcal{W}} c_t(x_{t+1}, u_t)\bigr],    
\end{flalign}
where the stage loss $c_t(x_{t+1}, u_t)$ is defined as follows.
\begin{flalign} \nonumber 
    c_t(x_{t+1}, u_t) :=& (x_{t+1}-d_{t+1})^T Q_{t+1}^i(x_{t+1}-d_{t+1})\\  & + \sum_{j = 1}^{N}\sum_{h = 1}^{N} (u_t^j)^T (R^i_t)_{jh}u^h_t,\label{eq:stage_loss}
\end{flalign}
Note that for decoupled dynamics, we add a random noise $w^i_t$ drawn from the distribution $\mathcal{W}^i$ with zero mean to the set of decoupled dynamics in equation~\eqref{eq:dyn_decoupled}. By defining $w_t := (w^1_t, \dots, w^N_t)$, the loss function of the LQ game with noisy decoupled dynamics can be written as the one in equation~\ref{eq:cost_lqg_report}.
Considering the loss function defined above, the definitions of Nash equilibrium, potential game, identical interest games for LQ games, (Definitions~\ref{def:ne},~\ref{def:potential}, and~\ref{def:identical_interest})  the game \emph{pseudo-gradient}, $\mathcal{G}(\gamma)$, and its Jacobian in equation~\eqref{eq:jac_j} do not change. Hence, Lemma~\ref{Lm:Lemma7} is still valid for the game with noisy dynamics. Note that the derivative of $c_t(x_{t+1}, u_t)$ with respect to $x_t$ or $u_t^j$ is similar to the game without noise. Therefore, we only investigate the derivative of $x_{t+1}$ and $u_t^i$ with respect to the decision variable $\gamma_{\tau}^j$ in the presence of the noise.\\
\textbf{Proposition~\ref{prop:coupled_i}  part (a)}: For the open-loop information structure, note that the decision-making variable is $\gamma_t^j = u_t^j$, and the game dynamics can be written as $x_{t+1} = a^tx_0 + \sum_{\tau=0}^{t}a^{t-\tau}[b^1u^1_{\tau} + b^2u^2_{\tau} + w_{\tau}] $. Thus, the derivative of $x_{t+1}$ and $u_t^i$ with respect to $u_{\tau}^j$ is as follows.
\begin{flalign*}
  \frac{\partial x_{t+1}}{\partial u_{\tau}^j} =&
  \begin{cases}
    A^{t-\tau}B^i &  t \ge \tau\\
    0 & t < \tau
  \end{cases} \\
  \frac{\partial u_{t}^i}{\partial u_{\tau}^j}  =&
  \begin{cases}
    1 &  t = \tau \ \text{and}  \ i=j \\
    0 & \text{otherwise}
  \end{cases}
\end{flalign*}

The derivatives above are similar to the ones without noise. This proves that the conditions in Proposition~\ref{prop:coupled_i}, part (a), are valid for the game with noise. \\
\textbf{Proposition~\ref{prop:coupled_i}  part (b)}: For the full-state feedback information structure, note that the decision-making variable is $\gamma_t^j = K_t^j$. The derivative of $x_{t+1}$ and $u_t^i$ with respect to $K_{\tau}^j$ can be derived using the chain rule and assuming the states are scalar (as mentioned in Example 3) as follows.
\begin{flalign*}
  \frac{\partial x_{t+1}}{\partial K_{\tau}^j} =&
  \begin{cases}
    -\prod_{t' = \tau+1}^{t} \frac{\partial x_{t'+1}}{\partial x_{t'}} b^jx_{\tau} & t > \tau\\
    - b^jx_t & t = \tau\\
    0 & t < \tau
  \end{cases} \\
  \frac{\partial u_{t}^i}{\partial K_{\tau}^j} =&
  \begin{cases}
  - K_t^i\frac{\partial x_{t}}{\partial K_{\tau}^j} & t> \tau \\
    - x_t &  t = \tau \ \text{and}  \ i=j \\
    0 & \text{otherwise}
  \end{cases}
\end{flalign*}
% \begin{flalign}
%     \frac{\partial x_{t+1}}{\partial K_{\tau}^j} &=-\prod_{t' = \tau+1}^{t} \frac{\partial x_{t'+1}}{\partial x_{t'}} b^jx_{\tau}  \qquad \text{If}  \  t > \tau; \qquad 
%     \frac{\partial x_{t+1}}{\partial K_{t}^j} = - b^jx_t  \qquad \  \text{If}  \ \tau = t; \qquad \text{else} \ \frac{\partial x_{t+1}}{ \partial K_{\tau}^j} = 0,\\ \frac{\partial u_{t}^i}{\partial K_{\tau}^j} &= - K_t^i\frac{\partial x_{t}}{\partial K_{\tau}^j}  \qquad \qquad \qquad  \text{If}  \ t > \tau; \qquad 
%     \frac{\partial u_{t}^i}{\partial K_{t}^j} = - x_t  \qquad \qquad  \text{If}  \ i = j; \qquad \text{else} \ \frac{\partial u_{t}^i}{\partial K_{\tau}^j} = 0 ,  
% \end{flalign}
Note that $\frac{\partial x_{t+1}}{\partial x_{t}}$ is not a function of $w_t$; therefore, the derivatives above are similar to the ones without noise. This proves that the conditions in Proposition 5, part (b), are valid for the game with noise.\\
\textbf{Theorem 9}: For the decoupled state feedback information structure, note that the decision-making variable is $\gamma_t^j = k_t^j$. The derivative of $x^i_{t+1}$ and $u_t^i$ with respect to $k_{\tau}^i$ can be derived using the chain rule as follows.
\begin{flalign*}
  \frac{\partial x^i_{t+1}}{\partial k_{\tau}^i} =&
  \begin{cases}
    -\prod_{t' = \tau+1}^{t} \frac{\partial x^i_{t'+1}}{\partial x^i_{t'}} b^i(x^i_{\tau})^T & t > \tau\\
    - b^i(x^i_t)^T & t = \tau\\
    0 & t < \tau
  \end{cases} \\
   \frac{\partial u_{t}^i}{\partial k_{\tau}^i} =&
  \begin{cases}
  - k_t^i\frac{\partial x^i_{t}}{\partial k_{\tau}^i}  & t> \tau \\
    - (x^i_t)^T &  t = \tau  \\
    0 & t < \tau
  \end{cases}
\end{flalign*}
% \begin{flalign}
%     \frac{\partial x^i_{t+1}}{\partial k_{\tau}^i} &=-\prod_{t' = \tau+1}^{t} \frac{\partial x^i_{t'+1}}{\partial x^i_{t'}} b^i(x^i_{\tau})^T  \qquad \text{If}  \  t > \tau;
%     \qquad 
%     \frac{\partial x^i_{t+1}}{\partial k_{t}^i} = - b^i(x^i_t)^T  \qquad \  \text{If}  \ \tau = t;
%     \qquad \text{else} \ \frac{\partial x^i_{t+1}}{ \partial k_{\tau}^i} = 0,\\ \frac{\partial u_{t}^i}{\partial k_{\tau}^i} &= - k_t^i\frac{\partial x^i_{t}}{\partial k_{\tau}^i}  \qquad \qquad \qquad \qquad \ \text{If}  \ t > \tau;  
%     \qquad 
%     \frac{\partial u^i_{t}}{\partial k_{t}^i} = - (x^i_t)^T  \qquad \qquad  \text{If}  \ \tau = t;
%     \qquad \text{else} \ \frac{\partial u_{t}^i}{\partial k_{\tau}^i} = 0. 
% \end{flalign}
Note that the derivative of $x^i_{t+1}$ and $u_t^i$ with respect to $k_t^j$ for $j \neq i$ are zero, and $\frac{\partial x^i_{t+1}}{\partial x^i_{t}}$ is not a function of $w_t$. Hence, the derivatives above are similar to the ones without noise. These results indicate that the conditions derived for the game with the decoupled dynamic
defined in equation~\eqref{eq:dyn_decoupled} and a decoupled state feedback information structure in (iii) to be potential do not depend on noise, and consequently, adding the noise does not affect the results.

\end{document}